\documentclass[10pt]{amsart}
\usepackage{mathptmx}
\usepackage{amsmath}     
\usepackage{amssymb}
\usepackage{array}
\usepackage{geometry}
\usepackage[bookmarks=true,colorlinks=true, pdfstartview=FitV, linkcolor=black, citecolor=blue, urlcolor=black]{hyperref}
\usepackage{movie15}

\usepackage{color}
\definecolor{DarkRed}{rgb}{0.55,.00,0.2}
\definecolor{DarkGrey}{rgb}{0.35,.35,0.35}

\theoremstyle{definition}

\theoremstyle{remark}

\numberwithin{equation}{section}

\begin{document}

\title{Orthogonal polynomials for the weight $x^{\nu} \exp(-x- t/x)$}

\author{{\bf S. Yakubovich} \\
 {\em {Department of Mathematics, Fac. Sciences of University of Porto,\\ Rua do Campo Alegre,  687; 4169-007 Porto (Portugal) }}}
  \thanks{ E-mail: syakubov@fc.up.pt} 
\vspace{3mm}

\thanks{ The work was  partially supported by CMUP, which is financed by national funds through FCT(Portugal),  under the project with reference UIDB/00144/2020. }


\date{\today}



\begin{abstract}  Orthogonal polynomials for the weight  $x^{\nu} \exp(-x- t/x),\ x, t > 0, \nu \in \mathbb{R}$ are investigated.  Differential-difference equations,  recurrence relations, explicit representations, generating functions and Rodrigues-type formula  are obtained.  
\end{abstract}

\maketitle

\markboth{\rm \centerline{ S. Yakubovich}}{orthogonal  polynomials for exponential weights}

\section{Introduction and preliminary results}

Let $\nu \in \mathbb{R},\  t > 0$ be parameters and consider the sequence of orthogonal polynomials $\{P_n^\nu(x,t)\}_{n\ge 0}$ of degree $n$, satisfying the orthogonality condition

$$\int_0^\infty  P_n^\nu(x,t) P_m^\nu(x,t) \ e^{-x-t/x} x^\nu dx=  \delta_{n,m},\eqno(1.1)$$ 
where $\delta_{n,m},\   n,m\in\mathbb{N}_{0}$ is the Kronecker symbol.  Admitting the limit case $t=0$ and $\nu > -1$, it gives normalized Laguerre polynomials

$$P_n^\nu(x,0) \equiv {\tilde L}_n^\nu (x) =  \left( {n!\over \Gamma(n+\nu+1)}\right)^{1/2} L_n^\nu (x),\eqno(1.2)$$
where $\{L_n^\nu (x)\}_{n\ge0}$ are classical Laguerre polynomials  \cite{Bateman}, Vol. II and $\Gamma(z)$ is the Euler gamma function  \cite{Bateman}, Vol. I.  It satisfies the three term recurrence relation

$$x L_n^\nu (x) = - (n+1) L_{n+1}^\nu (x)  + (2n+\nu+1) L_n^\nu (x) - (n+\nu) L_{n-1}^\nu (x),\  n \in \mathbb{N}.\eqno(1.3)$$
As it follows from the theory of orthogonal polynomials,  the three term recurrence relation for the sequence $\{P_n^\nu(x,t)\}_{n\ge 0}$ can be written in the form

$$x P_n^\nu(x,t) = A_{n+1} (t) P_{n+1}^\nu(x,t) +B_n(t) P_n^\nu(x,t) + A_{n}(t) P_{n-1}^\nu(x,t),\eqno(1.4)$$
where $P_{-1}^\nu(x,t) \equiv 0,\ P_n^\nu(x,t) = a_n(t) x^n+ b_n(t) x^{n-1}+ \dots,\quad a_n(t)\neq 0 $ and 

$$A_n(t)= {a_{n-1}(t)\over a_n(t)},\quad B_n(t)= {b_{n}(t)\over a_n(t)} - {b_{n+1}(t)\over a_{n+1}(t)}.\eqno(1.5)$$
As a consequence of (1.4) the Christoffel-Darboux formula takes place

$$\sum_{k=0}^n   P_k^\nu(x,t) P_k^\nu(y,t) = A_{n+1}(t) \frac{  P_{n+1}^\nu(x,t) P_n^\nu(y,t) -  P_{n}^\nu(x,t) P_{n+1}^\nu(y,t)}{x-y}.\eqno(1.6)$$
Using the value of the integral \cite{Bateman}, Vol. II

$$ \int_0^\infty  e^{-t/x -x}  x^{\nu-1} dx = 2 t^{\nu/2} K_\nu\left( 2\sqrt t\right) \equiv \rho_\nu(t),\ t > 0,\eqno(1.7)$$
where $K_\nu(z)$ is the modified Bessel function or Macdonald function \cite{YaL}, we can find the moments of the weight $x^{\nu} \exp(-x- t/x)$. Namely, we find

$$ \int_0^\infty  e^{-t/x -x}  x^{\nu+n} dx = \rho_{\nu+n+1}(t),\quad n \in \mathbb{N}_0.\eqno(1.8)$$
The asymptotic behavior of the modified Bessel function at infinity and near the origin \cite{Bateman}, Vol. II gives the corresponding values for the  function $\rho_\nu,\ \nu \in \mathbb{R}$.  Precisely, we have
$$\rho_\nu (t)= O\left( t^{(\nu-|\nu|)/2}\right),\  t \to 0,\ \nu\neq 0, \quad  \rho_0(t)= O( \log t),\ t \to 0,\eqno(1.9)$$

$$ \rho_\nu(t)= O\left( t^{\nu/2- 1/4} e^{- 2\sqrt t} \right),\ t \to +\infty.\eqno(1.10)$$
Moreover, it can be represented in terms of Laguerre polynomials (cf. \cite{PrudnikovMarichev}, Vol. II, Entry 2.19.4.13 )
$${(-1)^n t^n\over n!}\  \rho_\nu(t)=   \int_0^\infty x^{\nu+n -1} e^{-x - t/x}  L_n^\nu(x) dx,\quad    n \in\mathbb{N}_{0}.\eqno(1.11)$$
Further, it has a relationship with the  Riemann-Liouville fractional integral \cite{YaL}

$$ \left( I_{-}^\alpha  f \right) (t)  = {1\over \Gamma(\alpha)} \int_t^\infty (x-t)^{\alpha-1} f(x) dx,\quad  {\rm Re} \alpha > 0,\eqno(1.12)$$
namely,  we get  the  formula 
$$\rho_\nu(t)= \left( I_{-}^\nu \rho_0 \right) (t),\ \nu >0.\eqno(1.13)$$
Hence the index law for fractional integrals immediately implies

$$ \rho_{\nu+\mu} (t)= \left( I_{-}^\nu \rho_\mu \right) (t)=   \left( I_{-}^\mu \rho_\nu \right) (t).\eqno(1.14)$$
The corresponding definition of the fractional derivative presumes the relation $ D^\mu_{-}= - D  I_{-}^{1-\mu}$.   Hence for the ordinary $n$-th derivative of $\rho_\nu$ we find

$$D^n \rho_\nu(t)= (-1)^n \rho_{\nu-n} (t),\quad n \in \mathbb{N}_0.\eqno(1.15)$$
 Recalling (1.7) and integrating by parts, it is not difficult to establish  the following recurrence relation for $\rho_\nu$ 

$$\rho_{\nu+1} (t) =    \nu \rho_\nu(t)+ t \rho_{\nu-1} (t),\quad \nu \in \mathbb{R}.\eqno(1.16)$$
In the operator form it can be written as follows

$$\rho_{\nu+1} (t) = \left( \nu - tD \right) \rho_\nu(t).\eqno(1.17)$$
Finally in this section we observe that up to a normalization factor the orthogonality (1.1) is equivalent to the following $n$ equalities

$$\int_0^\infty  P_n^\nu(x,t) \ e^{-x-t/x} x^{\nu+m} dx=  0,\quad  m=0,1,\dots, n-1.\eqno(1.18)$$ 
Besides, taking into account (1.4),  (1.5), we get the identities   

$$\int_0^\infty  [P_n^\nu(x,t) ]^2 \ e^{-x-t/x} x^{\nu+1} dx= B_n(t),\eqno(1.19)$$

$$\int_0^\infty  P_n^\nu(x,t) \ e^{-x-t/x} x^{\nu+n} dx = {1\over a_n(t)},\eqno(1.20)$$

$$ \int_0^\infty  P_n^\nu(x,t) \ e^{-x-t/x} x^{\nu+n+1} dx = - { b_{n+1}(t)\over a_{n+1}(t) a_n(t)},\eqno(1.21)$$

$$ \int_0^\infty  P_n^\nu(x,t) P_{n-1}^\nu(x,t) \ e^{-x-t/x} x^{\nu+1} dx = A_n(t),\eqno(1.22)$$
and with (1.8) it yields

$$P_0^\nu(x,t) = [\rho_{\nu+1}(t)]^{-1/2},\quad   P_1^\nu(x,t) = \left[ \rho_{\nu+1}(t) (\rho_{\nu+3}(t)\rho_{\nu+1}(t) - \rho^2_{\nu+2}(t))\right]^{-1/2}$$

$$\times \left[ - \rho_{\nu+1}(t) x + \rho_{\nu+2}(t)\right],\eqno(1.23)$$ 
where $\rho_{\nu+3}(t)\rho_{\nu+1}(t) - \rho^2_{\nu+2}(t) > 0$ for all $t >0$. The latter fact follows immediately from the Ismail integral representation of the quotient of functions $\rho_\nu, \rho_{\nu+1}$ \cite{Ismail}

$${\rho_\nu(x) \over \rho_{\nu+1}(x) } = {1\over \pi^2} \int_0^\infty {y^{-1} dy \over (x+y) \left[ J_{\nu+1}^2 (2\sqrt y)+ Y^2_{\nu+1}(2\sqrt y) \right] }\eqno(1.24)$$
and the Nicholson's integral \cite{Bateman}, Vol. II 

$$J_{\nu}^2 (2\sqrt x)+ Y^2_{\nu}(2\sqrt x) = {8\over \pi^2} \int_0^\infty K_0\left( 4\sqrt x\ \sinh (t)\right) \cosh(2\nu t) dt.\eqno(1.25)$$

\section{Differential and differential-difference equations}

We begin this section with the following auxiliary lemma.

{\bf Lemma 1}. {\it Let $\nu \in \mathbb{R},\  t > 0,\ n \in \mathbb{N}$.   Then it has the identities}

$$\int_0^\infty \left[P_n^\nu(x,t)\right]^2  \ e^{-x-t/x} x^{\nu-1} dx =  {1\over t} \left[ B_n(t)-\nu-1-2n \right],\eqno(2.1)$$

$$\int_0^\infty \left[P_n^\nu(x,t)\right]^2  \ e^{-x-t/x} x^{\nu-2} dx = {1\over t} \left[ 1- {\nu\over t} \left[ B_n(t)-\nu-1-2n \right] \right],\eqno(2.2)$$

$$\int_0^\infty P_n^\nu(x,t) P_{n-1}^\nu(x,t)  \ e^{-x-t/x} x^{\nu-1} dx = {1\over t} \left[ A_n(t) + { b_{n}(t)\over a_{n}(t) A_n(t)} \right],\eqno(2.3)$$

$$\int_0^\infty P_n^\nu(x,t) P_{n-1}^\nu(x,t)  \ e^{-x-t/x} x^{\nu-2} dx = - {1\over t} \left[  {\nu\over t} \left[ A_n(t) + { b_{n}(t)\over a_{n}(t) A_n(t)} \right] + {n\over A_n(t)} \right].\eqno(2.4)$$

\begin{proof} In fact, employing (1.18),  (1.19), (1.20) and integration by parts we derive

$$\int_0^\infty \left[P_n^\nu(x,t)\right]^2  \ e^{-x-t/x} x^{\nu-1} dx = {1\over t} \left[  \int_0^\infty \left[P_n^\nu(x,t)\right]^2  \ e^{-x-t/x} x^{\nu+1} dx\right.$$

$$\left. - (\nu+1) \int_0^\infty \left[P_n^\nu(x,t)\right]^2  \ e^{-x-t/x} x^{\nu} dx - 2 \int_0^\infty P_n^\nu(x,t) {\partial\over \partial x} [P_n^\nu(x,t)]  \ e^{-x-t/x} x^{\nu+1} dx\right] $$

$$= {1\over t} \left[ B_n(t)-\nu-1-2n \right],$$
which proves (2.1).  Then, analogously, recalling (1.1) and using (2.1), we get

$$\int_0^\infty \left[P_n^\nu(x,t)\right]^2  \ e^{-x-t/x} x^{\nu-2} dx = {1\over t} \left[  \int_0^\infty \left[P_n^\nu(x,t)\right]^2  \ e^{-x-t/x} x^{\nu} dx\right.$$

$$\left.  - \nu \int_0^\infty \left[P_n^\nu(x,t)\right]^2  \ e^{-x-t/x} x^{\nu-1} dx - 2 \int_0^\infty P_n^\nu(x,t) {\partial\over \partial x} [P_n^\nu(x,t)]  \ e^{-x-t/x} x^{\nu} dx\right]$$

$$ = {1\over t} \left[ 1- {\nu\over t} \left[ B_n(t)-\nu-1-2n \right] \right].$$
Concerning (2.3), we deduce, invoking (1.1), (1.5), (1.20), (1.21), (1.22),   

$$\int_0^\infty P_n^\nu(x,t) P_{n-1}^\nu(x,t)  \ e^{-x-t/x} x^{\nu-1} dx = {1\over t} \left[  \int_0^\infty P_n^\nu(x,t) P_{n-1}^\nu(x,t)   \ e^{-x-t/x} x^{\nu+1} dx \right.$$

$$\left. - (\nu+1) \int_0^\infty P_n^\nu(x,t)  \ P_{n-1}^\nu(x,t)  e^{-x-t/x} x^{\nu} dx - \int_0^\infty P_n^\nu(x,t) {\partial\over \partial x} [P_{n-1}^\nu(x,t)]  \ e^{-x-t/x} x^{\nu+1} dx\right.$$

$$\left. -  \int_0^\infty P_{n-1}^\nu(x,t) {\partial\over \partial x} [P_n^\nu(x,t)]  \ e^{-x-t/x} x^{\nu+1} dx\right] = {1\over t} \left[ A_n(t) + { b_{n}(t)\over a_{n}(t) A_n(t)} \right].$$
Finally,  a similar reasoning gives via (2.3)

$$\int_0^\infty P_n^\nu(x,t) P_{n-1}^\nu(x,t)  \ e^{-x-t/x} x^{\nu-2} dx = {1\over t} \left[  \int_0^\infty P_n^\nu(x,t) P_{n-1}^\nu(x,t)   \ e^{-x-t/x} x^{\nu} dx \right.$$

$$\left. - \nu \int_0^\infty P_n^\nu(x,t)  \ P_{n-1}^\nu(x,t)  e^{-x-t/x} x^{\nu-1} dx - \int_0^\infty P_n^\nu(x,t) {\partial\over \partial x} [P_{n-1}^\nu(x,t)]  \ e^{-x-t/x} x^{\nu} dx\right.$$

$$\left. -  \int_0^\infty P_{n-1}^\nu(x,t) {\partial\over \partial x} [P_n^\nu(x,t)]  \ e^{-x-t/x} x^{\nu} dx\right] = - {1\over t} \left[  {\nu\over t} \left[ A_n(t) + { b_{n}(t)\over a_{n}(t) A_n(t)} \right] + {n\over A_n(t)} \right],$$
and we establish (2.4), completing the proof of Lemma 1. 

\end{proof}

Now we are ready to prove the following theorem.

{\bf Theorem 1}. {\it Let $n \in \mathbb{N}$. Orthogonal polynomials $P_n^\nu(x,t)$ satisfy the first order linear differential-difference equation }

$$x^2 {\partial\over \partial x} [P_n^\nu(x,t)]  = \left[ nx - A_n^2(t) - {b_n(t)\over a_n(t)} \right] P_n^\nu(x,t) + A_n(t) \left[ x+B_n(t)-\nu-1-2n\right] P_{n-1}^\nu(x,t).\eqno(2.5)$$

\begin{proof} Since $ {\partial\over \partial x} [P_n^\nu(x,t)] $ is a polynomial of degree $n-1$, we write it in the form

$${\partial\over \partial x} [P_n^\nu(x,t)]  = \sum_{k=0}^{n-1} c_{n,k} (t) P_k^\nu(x,t),\eqno(2.6)$$
where, owing to the orthogonality,

$$c_{n,k} (t)= \int_0^\infty  {\partial\over \partial x} [P_n^\nu(x,t)] P_k^\nu(x,t) \ e^{-x-t/x} x^\nu dx.\eqno(2.7)$$
Then, integrating by parts and using the orthogonality relation, we obtain

$$c_{n,k} (t)= - t  \int_0^\infty  P_n^\nu(x,t) P_k^\nu(x,t) \ e^{-x-t/x} x^{\nu-2} dx - \nu  \int_0^\infty  P_n^\nu(x,t) P_k^\nu(x,t) \ e^{-x-t/x} x^{\nu-1} dx.\eqno(2.8)$$
Moreover, we observe that 

$$ \int_0^\infty  P_n^\nu(y,t) \sum_{k=0}^{n-1}  P_k^\nu(y,t) P_k^\nu(x,t) \ e^{-y-t/y} y^{\nu} x^{-1} dy = 0,$$

$$ \int_0^\infty  P_n^\nu(y,t) \sum_{k=0}^{n-1}  P_k^\nu(y,t) P_k^\nu(x,t) \ e^{-y-t/y} y^{\nu} x^{-2} dy = 0.$$
Therefore from (2.6), (2.8) and Christoffel-Darboux formula (1.6) we derive

$${\partial\over \partial x} [P_n^\nu(x,t)]  = -t  \sum_{k=0}^{n-1} P_k^\nu(x,t) \int_0^\infty  P_n^\nu(y,t) P_k^\nu(y,t) \ e^{-y-t/y} y^{\nu} \left[ {1\over y^2} - {1\over x^2} \right] dy$$

$$- \nu \sum_{k=0}^{n-1} P_k^\nu(x,t) \int_0^\infty  P_n^\nu(y,t) P_k^\nu(y,t) \ e^{-y-t/y} y^{\nu} \left[ {1\over y} - {1\over x} \right] dy$$

$$= -t  A_{n}(t) P_{n}^\nu(x,t) \int_0^\infty  P_n^\nu(y,t) P_{n-1}^\nu(y,t) \ e^{-y-t/y} y^{\nu} \left[ {1\over x y^2} + {1\over y x^2} \right] dy$$

$$+ t  A_{n}(t) P_{n-1}^\nu(x,t) \int_0^\infty  \left[P_{n}^\nu(y,t)\right]^2  \ e^{-y-t/y} y^{\nu} \left[ {1\over x y^2} + {1\over y x^2} \right] dy$$

$$ - {\nu\over x}   A_{n}(t) P_{n}^\nu(x,t) \int_0^\infty  P_n^\nu(y,t) P_{n-1}^\nu(y,t) \ e^{-y-t/y} y^{\nu-1} dy$$

$$+ {\nu\over x}   A_{n}(t) P_{n-1}^\nu(x,t) \int_0^\infty  \left[P_{n}^\nu(y,t)\right]^2  \ e^{-y-t/y} y^{\nu-1} dy.\eqno(2.9)$$
Hence by virtue of Lemma 1 equalities (2.9) become

$${1\over A_n(t)} {\partial\over \partial x} [P_n^\nu(x,t)]  =  {1\over x}   P_{n}^\nu(x,t) \left[  \left( {\nu\over t} - {1\over x} \right) \left[ A_n(t) + { b_{n}(t)\over a_{n}(t) A_n(t)} \right] + {n\over A_n(t)}   \right]$$

$$ + {1\over x}   P_{n-1}^\nu(x,t) \left[ 1- \left( {\nu\over t} - {1\over x}\right) \left[ B_n(t)-\nu-1-2n \right]  \right]$$

$$- {\nu \over  x t} P_{n}^\nu(x,t) \left[ A_n(t) + { b_{n}(t)\over a_{n}(t) A_n(t)} \right] + {\nu\over x t}  P_{n-1}^\nu(x,t) \left[ B_n(t)-\nu-1-2n \right].$$
Hence after simplification we arrive at the differential-difference equation (2.5).  Theorem 1 is proved. 
\end{proof} 

{\bf Corollary 1.}  {\it  Denoting by $a_{n,0}(t)$ the free term of the polynomial $P_n^\nu(x,t)$,  it has the value  }

$$a_{n,0}(t)= {1\over a_n(t) \rho_{\nu+1}(t) } \prod_{k=1}^n \frac{B_k(t) -\nu-1-2k} {A^2_k(t)+ {b_k(t)\over a_k(t)}}.\eqno(2.10)$$

\begin{proof} In fact, letting $x=0$ in (2.5), we find the recurrence relation

$$\left[A_n^2(t) + {b_n(t)\over a_n(t)} \right] a_{n,0}(t) = \left[ B_n(t)-\nu-1-2n\right] A_n(t) a_{n-1,0}(t).\eqno(2.11)$$
Hence formula (2.10) comes immediately, solving  recurrence (2.11) with the use of (1.5) and (1.23). 

\end{proof}

{\bf Remark 1}.  In the limit case $t=0$ we use (1.2), (1.5), (1.8) to have  $\rho_{\nu+1}(0) = \Gamma(\nu+1)$,

$$a_n(0) = {(-1)^n \over[ n! \Gamma(n+\nu+1)]^{1/2}},\quad b_n(0) = (-1)^{n+1}  \left( {n(n+\nu) \over (n-1)! \Gamma(n+\nu)}\right)^{1/2},\eqno(2.12) $$ 

$$B_n(0)= 2n+\nu+1,\quad  A_n(0) = - (n(n+\nu))^{1/2},\quad a_{n,0}(0) = {1\over \Gamma(\nu+1) }\left( { \Gamma(n+\nu+1)\over n!}\right)^{1/2},$$

$$\lim_{t\to 0}  \left[A_n^2(t) + {b_n(t)\over a_n(t)} \right] \bigg[ B_n(t)-\nu-1-2n\bigg]^{-1} = -n.\eqno(2.13)$$
Moreover,   taking the limit  of the indeterminate form under the product sign in (2.10) when $t\to 0$, it implies the identity 

$$\prod_{k=1}^n \frac{B^\prime_k(0) } { (-1)^k (k! \Gamma(k+\nu+1))^{1/2} [b^\prime_k(0) + k(k+\nu) a^\prime_k(0)] - 2  ( k (k+\nu))^{1/2} A^\prime_k(0)} =  {(-1)^n\over n!}.\eqno(2.14)$$

Now, assuming that polynomial coefficients in (1.4), (1.5) $a_{n,k}(t) \in  C^1(\mathbb{R}_+),\ a_{n,n}(t) \equiv a_n(t),\ a_{n,n-1}(t)= b_n(t)$ as functions of $t$, we differentiate equality (1.1) by $t$ to obtain  

$$0= {\partial\over \partial t} \int_0^\infty \left[P_n^\nu(x,t)\right]^2  \ e^{-x-t/x} x^{\nu} dx = 2 \int_0^\infty  P_n^\nu(x,t) {\partial\over \partial t} \left[P_n^\nu(x,t)\right]  \ e^{-x-t/x} x^{\nu} dx $$

$$- \int_0^\infty \left[P_n^\nu(x,t)\right]^2  \ e^{-x-t/x} x^{\nu-1} dx,$$
where the differentiation under the integral sign can be easily motivated by virtue of the absolute and uniform convergence.  Then, appealing to (1.20) and (2.1), it implies the equality

$$ {a^\prime_n(t)\over a_n(t)} = {B_n(t)-\nu-1-2n\over 2t}.\eqno(2.15)$$ 
This first order differential equation with respect to the leading term $a_n(t)$ can be uniquely solved under  initial condition (2.12), and we get

$$a_n(t)=    {(-1)^n \over \left(n! \Gamma(n+\nu+1)\right)^{1/2}}\ \exp\left({1\over 2} \int_0^t {B_n(y)-\nu-1-2n\over y} dy\right),\eqno(2.16)$$
where the integral under the exponential function exists since, evidently, (see (2.13), (2.15))

$$\lim_{t\to 0} {B_n(t)-\nu-1-2n\over t} = B_n^\prime(0) = 2 {a^\prime_n(0)\over a_n(0)} .\eqno(2.17)$$
Moreover, employing (1.5),  we find

$$A_{n+1}(t)=  - ((n+1)(n+1+\nu))^{1/2} \ \exp\left({1\over 2} \int_0^t {B_n(y)- B_{n+1}(y) +2\over y} dy\right).\eqno(2.18)$$
Hence after differentiation we obtain 

$$A^\prime_{n+1}(0)=  {1\over 2} ((n+1)(n+1+\nu))^{1/2} \left[ B^\prime_{n+1}(0)- B^\prime_n(0) \right].\eqno(2.19)$$
On the other hand,  via (1.15), (1.23) we have 

$$B_0(t)= {\rho_{\nu+2}(t)\over \rho_{\nu+1}(t)},\quad B^\prime_0(t)=  {\rho_{\nu+2}(t)\rho_\nu(t)\over \rho^2_{\nu+1}(t)}-1,\quad  B^\prime_0(0) = {1\over \nu},\ \nu > 0.\eqno(2.20)$$
Hence it gives the equality 
 
 $$B^\prime_{n+1}(0) = {1\over \nu}+ 2 \sum_{k=0}^n {A^\prime_{k+1}(0)\over ((k+1)(k+1+\nu))^{1/2}} .\eqno(2.21)$$

The following theorem shows the first order partial differential-difference equation whose solutions are polynomials $P_n^\nu(x,t)$.

{\bf Theorem 2}. {\it Orthogonal polynomials $P_n^\nu(x,t)$ satisfy the first order partial  differential-difference equation }

$$\left( t {\partial\over \partial t} + x {\partial\over \partial x} \right) P_n^\nu(x,t)  = \left( t {a_n^\prime(t)\over a_n(t)}+ n \right) P_n^\nu(x,t) + A_n(t) P_{n-1}^\nu(x,t).\eqno(2.22)$$
\begin{proof} Indeed, differentiating both sides of (1.18) with respect to $t$, we have

$$\int_0^\infty  {\partial\over \partial t} [ P_n^\nu(x,t) ] \ e^{-x-t/x} x^{\nu+m} dx -  \int_0^\infty   P_n^\nu(x,t) \ e^{-x-t/x} x^{\nu+m-1} dx =  0,\quad  m=0,1,\dots, n-1.\eqno(2.23)$$ 
Meanwhile, the second integral on the left-hand side of (2.23) can be treated via integration by parts which gives 

$$ \int_0^\infty   P_n^\nu(x,t) \ e^{-x-t/x} x^{\nu+m-1} dx =  {1\over t} \int_0^\infty   P_n^\nu(x,t) \ e^{-x-t/x} x^{\nu+m+1} dx$$

$$ -  {1\over t} \int_0^\infty  {\partial\over \partial x} [ P_n^\nu(x,t) ] \ e^{-x-t/x} x^{\nu+m+1} dx  - {\nu+m+1\over t}  \int_0^\infty   P_n^\nu(x,t) \ e^{-x-t/x} x^{\nu+m} dx,$$
and, combining with (2.23), (1.4) and (1.18), we obtain  

$$ \int_0^\infty  \left[ \left( t {\partial\over \partial t} + x {\partial\over \partial x} \right) P_n^\nu(x,t) - A_n(t) P_{n-1}^\nu(x,t) \right] \ e^{-x-t/x} x^{\nu+m} dx =  0,\quad  m=0,1,\dots, n-1.\eqno(2.24)$$ 
Hence by unicity we therefore have

$$ \left( t {\partial\over \partial t} + x {\partial\over \partial x} \right) P_n^\nu(x,t) - A_n(t) P_{n-1}^\nu(x,t) = c_n(t) P_n^\nu(x,t).\eqno(2.25)$$
The function $c_n(t)$ is defined, equating  coefficients of $x^n$ on both sides in (2.25) to find

$$c_n(t)= t {a^\prime_n(t)\over a_n(t)} + n.$$
This completes the proof of Theorem 2.

\end{proof} 

{\bf Corollary 2}. {\it Equation $(2.22)$ can be written in the form}

$$\left( t {\partial\over \partial t} + x {\partial\over \partial x} \right) P_n^\nu(x,t)  = {1\over 2} \left[B_n(t)-\nu-1\right] P_n^\nu(x,t) + A_n(t) P_{n-1}^\nu(x,t).$$

\begin{proof} The proof is immediate with the use of (2.15).

\end{proof}

{\bf Corollary 3}. {\it Let $t >0$. The following equalities take place}

$${d\over dt} \left[{b_n(t)\over a_n(t)}\right] =  {1\over t} \left[ A_n^2(t) + {b_n(t)\over a_n(t)}\right],\eqno(2.26)$$

$$B_n^\prime(t)= {1\over t} \left[ A_n^2(t) - A_{n+1}^2(t)+ B_n(t)\right],\eqno(2.27)$$

$$a_{n,0}^\prime(t)=  {B_n(t)-\nu-1\over 2t} a_{n,0}(t) +  {A_n(t)\over t}  a_{n-1,0}(t).\eqno(2.28)$$

\begin{proof} The proof is immediate, equating coefficients of $x^{n-1}$ and free terms on both sides of (2.22), and the use of (1.5), (2.15).

\end{proof} 
Further, differentiating (1.19) by $t$ and employing (1.1), (1.21), we deduce

$$B_n^\prime(t)= 2 \int_0^\infty P_n^\nu(x,t) {\partial\over \partial t} [ P_n^\nu(x,t) ] \ e^{-x-t/x} x^{\nu+1} dx - 1$$

$$= - 2 {a_n^\prime(t) b_{n+1}(t)\over a_{n+1}(t) a_n(t)} + 2 {b_n^\prime(t) \over a_n(t)} -1,$$
i.e. (see (1.5))

$$ B_n^\prime(t)= 2 B_n(t)  {a_n^\prime(t) \over a_n(t)} + 2 {d\over dt} \left[{b_n(t)\over a_n(t)}\right] -1.\eqno(2.29)$$

{\bf Corollary 4}. {\it Let $\nu > 0,\   n \in \mathbb{N}_0$. It has the values }

$$A_n^\prime(0) = 0,\quad  B_n^\prime (0)= {1\over \nu},\quad a_n^\prime(0)=  {(-1)^n \over 2\nu [n! \Gamma(n+\nu+1)]^{1/2}},\eqno(2.30)$$

$$    b_n^\prime(0)=  (-1)^{n+1} {n(n+\nu+2)  \over 2\nu [n! \Gamma(n+\nu+1)]^{1/2}},\quad  a_{n,0}^\prime(0)= {1\over 2\nu \Gamma(\nu+1) }\left( { \Gamma(n+\nu+1)\over n!}\right)^{1/2}.\eqno(2.31)$$

\begin{proof} In fact, taking a limit in (2.26) when $t\to 0$, we find

$$\lim_{t\to 0} {d\over dt} \left[ A_n^2(t) \right] = 2 A_n(0) A^\prime_n(0) = 0.$$
Hence via (2.13) $A_n^\prime(0) = 0$ and the sum in (2.21) is zero. Thus  $B_n^\prime (0)= {1\over \nu}$. The value for $a_n^\prime(0)$ comes from (2.12), (2.17) and  $b_n^\prime (0) $ is obtained from  (2.29). Finally, $a_{n,0}^\prime(0)$ is a consequence of (2.28) when $t \to 0$.

\end{proof} 
Now, assuming that polynomial coefficients are twice continuously differentiable functions of $t$, we solve a simple Cauchy problem for the first order differential equation (2.26) to find

$$  {b_n(t)\over a_n(t)} = 2t \int_0^t {A_n(y) A_n^\prime(y)\over y} dy -  A^2_n(t) - {n t\over \nu}.\eqno(2.32)$$
Moreover, via (1.5), (2.27) and integration by parts we derive

$$  B_n(t)=   t \int_0^t {\left(A^2_n(y) - A^2_{n+1}(y)\right)^\prime\over y} dy + A^2_{n+1}(t) -  A^2_n(t) + { t\over \nu}$$

$$= 2n+\nu+1 +  t\int_0^t \left[A^2_n(y) - A^2_{n+1}(y)+2n+\nu+1\right] {dy \over y^2} + { t\over \nu},\eqno(2.33)$$

$$B^\prime_n(t)=   \int_0^t {\left(A^2_n(y) - A^2_{n+1}(y)\right)^\prime\over y} dy + {1\over \nu}.\eqno(2.34)$$

{\bf Corollary 5}. {\it Coefficients $A_n$ satisfy the following second kind nonlinear differential-difference equation}

$$   A_n^{\prime\prime}(t)  A_n(t) - \left(A_n^{\prime}(t)\right)^2  - {A^2_n(t)\over 2t^2}  \bigg[ A_{n-1}^2(t) - 2 A^2_n(t) + A_{n+1}^2(t) -2\bigg] = 0.\eqno(2.35)$$

\begin{proof} Indeed, from (1.5) and (2.18) we get

$${ A_n^\prime(t)\over A_n(t)} =  {B_{n-1}(t)- B_{n}(t) +2\over 2t},\quad t >0.\eqno(2.36)$$
Then by virtue of (2.27)

$$ 2 t {d\over dt} \left[ t { A_n^\prime(t)\over A_n(t)} \right] = t {d\over dt} \left[ B_{n-1}(t)- B_{n}(t)\right] =  B_{n-1}(t)- B_{n}(t) + A_{n-1}^2(t) - 2 A^2_n(t) + A_{n+1}^2(t)$$
and the result follows via (2.36) and simple differentiation.

\end{proof}

{\bf Corollary 6}. {\it The free term $a_{n,0}(t)$ can be determined by the formula

$$a_{n,0}(t) = (-1)^n (1+\nu)_n\ a_n(t) \exp\left( \int_0^t \left[ n  +  {A_n^2(y)  + {b_n(y)\over a_n(y)}\over  B_n(y)-\nu-1-2n } \right] {dy\over y} \right),\eqno(2.37)$$
where $(z)_n$ is the Pochhammer symbol}.

\begin{proof} In fact, invoking (2.11), (2.15), (2.28) and since $ B_n(t)-\nu-1-2n \neq 0,\ t >0$ via (2.1), we write the differential equation 

$${ a_{n,0}^\prime(t)\over a_{n,0}(t)} =  { a_{n}^\prime(t)\over a_{n}(t)} + {1\over t} \left[ n  +  {A_n^2(t)  + {b_n(t)\over a_n(t)}\over  B_n(t)-\nu-1-2n } \right],$$
which can be uniquely solved by formula (2.37), owing to inicial conditions (2.12) and where the integral converges under condition (2.13).

\end{proof}

{\bf Theorem 3}.   {\it Orthogonal polynomials $P_n^\nu(x,t)$ obey  the second order differential equation}

$$ x^4 \bigg[ x+B_n(t)-\nu-1-2n \bigg]  {\partial^2\over \partial x^2} [P_n^\nu(x,t)]  - x^2 \bigg[ x^2 + \left[ x+B_n(t)-\nu-1-2n \right] \bigg.$$

$$ \times \left[ (2n-3) x-  A_n^2(t) - A_{n-1}^2(t)- {b_n(t)\over a_n(t)}-  {b_{n-1}(t)\over a_{n-1} (t)}\right.$$

$$\bigg.\bigg.+  \left[ x+B_{n-1}(t)-\nu+1-2n\right] [x- B_{n-1}(t)]\bigg]\bigg]  {\partial\over \partial x} [ P_n^\nu(x,t)] $$

$$+ \bigg[    \left[ x+B_n(t)-\nu-1-2n\right] \bigg[  A^2_n(t)  \left[ x+B_n(t)-\nu-1-2n\right] \left[ x+B_{n-1}(t)-\nu+1-2n\right]\bigg.\bigg.$$

$$\bigg. + \left[nx-   A_n^2(t) - {b_n(t)\over a_n(t)} \right] \bigg[ (n-1)x-  A_{n-1}^2(t) - {b_{n-1}(t)\over a_{n-1}(t)}  + \left[ x+B_{n-1}(t)-\nu+1-2n\right] [x- B_{n-1}(t)]  \bigg]  \bigg]$$

$$\bigg. \bigg. - x^2 \bigg[ n \left[B_n(t)-\nu-2n -1\right]+  A_n^2(t) + {b_n(t)\over a_n(t)}  \bigg] \bigg]P_n^\nu(x,t) = 0.\eqno(2.38)$$

\begin{proof} Differentiating both sides of (2.5) with respect to $x$, we have 

$$x^2 {\partial^2\over \partial x^2} [P_n^\nu(x,t)]  = n P_n^\nu(x,t) +  \left[ (n-2) x - A_n^2(t) - {b_n(t)\over a_n(t)} \right] {\partial\over \partial x} [ P_n^\nu(x,t)]  $$

$$ + A_n(t) P_{n-1}^\nu(x,t)+ A_n(t) \left[ x+B_n(t)-\nu-1-2n\right] {\partial\over \partial x} [ P_{n-1}^\nu(x,t)].$$
But,  owing to (1.4) and (2.5), the latter equality becomes

$$x^4 {\partial^2\over \partial x^2} [P_n^\nu(x,t)]  =  x^2\left[ (n-2) x - A_n^2(t) - {b_n(t)\over a_n(t)} \right] {\partial\over \partial x} [ P_n^\nu(x,t)]  $$

$$+  \left[ n x^2 -   A^2_n(t) \left[ x+B_n(t)-\nu-1-2n\right] \left[ x+B_{n-1}(t)-\nu+1-2n\right] \right] P_{n}^\nu(x,t) $$

$$ + \left[ {x^2\over  x+B_n(t)-\nu-1-2n} + (n-1)x - A_{n-1}^2(t) - {b_{n-1}(t)\over a_{n-1}(t)} \right.$$

$$\bigg.  + \left[ x+B_{n-1}(t)-\nu+1-2n\right] [x- B_{n-1}(t)] \bigg] A_n(t) \left[ x+B_n(t)-\nu-1-2n\right]  P_{n-1}^\nu(x,t).$$
Finally, recalling (2.5) to express $P_{n-1}^\nu(x,t)$, we end up with the equation (2.37), completing the proof of Theorem 3.

\end{proof}

\section{Explicit representations. Recurrence relations for coefficients}

In this section we will deduce recurrence relations for  coefficients of  orthogonal polynomials $P_n^\nu(x,t)$ and their explicit representations.   In fact, we have

{\bf Theorem 4}.    {\it Let $n \in \mathbb{N}_0,\  t >0$. Then the following identities hold}

$$A_{n+1}^2(t)+ B_n^2(t)+ A_n^2(t)- (2n+\nu+2)B_n(t) + 2{b_n(t)\over a_n(t)} - t=0,\eqno(3.1)$$

$$A^2_n(t)  \left[B_n(t)-\nu-1-2n\right] \left[ B_{n-1}(t)-\nu+1-2n\right] - \left[   A_n^2(t) + {b_n(t)\over a_n(t)} \right] ^2+ t \bigg[    A_n^2(t) + {b_n(t)\over a_n(t)}  \bigg]  = 0,\eqno(3.2)$$

$$t \bigg[ B_{n-1}(t) + B_n(t) \bigg]^\prime + \bigg[ B_{n-1}(t) + B_n(t) \bigg] \bigg[ B_{n-1}(t) - B_n(t) +1\bigg] - (2n+\nu) \bigg[ B_{n-1}(t) - B_n(t) \bigg] =0.\eqno(3.3)$$

\begin{proof} In order to prove (3.1), we appeal to the three term recurrence relation (1.4) and the orthogonality condition (1.1). Hence, integrating by parts, we  derive

$$A_{n+1}^2(t)+ B_n^2(t)+ A_n^2(t)= \int_0^\infty  [x P_n^\nu(x,t) ]^2 \ e^{-x-t/x} x^{\nu} dx= (\nu+2) \int_0^\infty  [ P_n^\nu(x,t) ]^2 \ e^{-x-t/x} x^{\nu+1} dx $$

$$+ t \int_0^\infty  [P_n^\nu(x,t) ]^2 \ e^{-x-t/x} x^{\nu} dx + 2 \int_0^\infty   P_n^\nu(x,t) {d\over dx} [ P_n^\nu(x,t) ] \ e^{-x-t/x} x^{\nu+2} dx $$

$$= (\nu+2) B_n(t) + t + 2n B_n(t) + 2 A_n(t) \int_0^\infty  x  P_{n-1}^\nu(x,t) {d\over dx} [ P_n^\nu(x,t) ] \ e^{-x-t/x} x^{\nu} dx $$

$$t+ (2n+\nu+2) B_{n}(t) -  2n\ A_{n}(t) { b_n(t)\over  a_{n-1} (t)} + 2 (n-1) A_n(t)\ {b_n(t)\over a_{n-1}(t)} $$

$$= t+ (2n+\nu+2) B_{n}(t) - 2 {b_n(t)\over a_{n}(t)}.$$
This gives (3.1). On the other hand, writing (3.1) for the index $n-1$ in the form 

$$ A_n^2(t) + {b_n(t)\over a_n(t)}  -t =   B_{n-1}(t) \left[ 2n+\nu-1 - B_{n-1}(t)\right] -  A_{n-1}^2(t) - {b_{n-1}(t)\over a_{n-1}(t)},\eqno(3.4)$$
we let $x= 0$ in (2.38) to find the equality 

$$A^2_n(t)  \left[ B_n(t)-\nu-1-2n\right] \left[ B_{n-1}(t)-\nu+1-2n\right]\bigg.\bigg.$$

$$\bigg. + \left[  A_n^2(t) + {b_n(t)\over a_n(t)} \right] \bigg[   A_{n-1}^2(t) + {b_{n-1}(t)\over a_{n-1}(t)}  +\left[  B_{n-1}(t)-2n-\nu+1\right]  B_{n-1}(t)  \bigg]  = 0.\eqno(3.5)$$
Hence a simple comparison leads to (3.2). Finally, writing (3.1) for $n-1$ and subtracting one equality from another,  we use (2.27) to establish (3.3). 

\end{proof}

{\bf Corollary 7}.   {\it Differential equation $(2.38)$ reduces to the equality }

$$ x^2 \left( x+B_n(t)-\nu-1-2n \right)  {\partial^2\over \partial x^2} [P_n^\nu(x,t)] $$

$$ -  \bigg[x^3 +   \bigg[ B_n(t)- 2(\nu+n+1) \bigg] x^2 - \bigg[ t+ \left(B_n(t)-\nu-1-2n\right)(\nu+2) \bigg] x -  t \left(B_n(t)-\nu-1-2n\right)\bigg]  {\partial\over \partial x} [ P_n^\nu(x,t)] $$

$$+  \bigg[  n x^2 -   \left( {b_n(t)\over a_n(t)} - n \left( B_n(t)- 2\nu-3n-1\right) \right)  x +  \left( B_n(t)-\nu-2n-1\right) \bigg[ A^2_n(t)  - n (n+\nu+1) \bigg] \bigg.$$

$$\bigg. +  \left[  A_n^2(t) +{b_n(t)\over a_n(t)} \right] \left(2n+\nu - B_n(t)\right) \bigg] P_n^\nu(x,t) = 0.\eqno(3.6)$$

\begin{proof}  In fact, appealing to (3.4), (3.5), we get from (2.38)

$$ x^3 \bigg[ x+B_n(t)-\nu-1-2n \bigg]  {\partial^2\over \partial x^2} [P_n^\nu(x,t)]  - x \bigg[ x^2 + \left[ x+B_n(t)-\nu-1-2n \right] \bigg.$$

$$ \times \left[ (2n-3) x-  A_n^2(t) - A_{n-1}^2(t)- {b_n(t)\over a_n(t)}-  {b_{n-1}(t)\over a_{n-1} (t)}\right.$$

$$\bigg.\bigg.+  \left[ x+B_{n-1}(t)-\nu+1-2n\right] [x- B_{n-1}(t)]\bigg]\bigg]  {\partial\over \partial x} [ P_n^\nu(x,t)] $$

$$+ \bigg[    \left[ x+B_n(t)-\nu-1-2n\right] \bigg[  A^2_n(t)  \left[ B_n(t)+B_{n-1}(t)+x  -2 \nu- 4n\right] +  n( x^2 - (n+\nu) x -t ) \bigg] $$

$$ + \left[  A_n^2(t) +{b_n(t)\over a_n(t)} \right] \bigg[ \left[ B_n(t)-\nu-1-2n\right] (2n+\nu-x) + x(2n+\nu-1 - x) \bigg] $$

$$\bigg. - n x \left[B_n(t)-\nu-2n -1\right] \bigg]P_n^\nu(x,t) = 0.\eqno(3.7)$$
Hence, letting $x=0$ in (3.7), we find the identity

$$A^2_n(t)  \left[ B_n(t)+B_{n-1}(t) -2 \nu- 4n\right] + (2n+\nu) \left[  A_n^2(t) +{b_n(t)\over a_n(t)} \right] -nt =0.\eqno(3.8)$$
Consequently,  (3.7) becomes

$$ x^3 \bigg[ x+B_n(t)-\nu-1-2n \bigg]  {\partial^2\over \partial x^2} [P_n^\nu(x,t)]  -  x\bigg[ x^2 + \left[ x+B_n(t)-\nu-1-2n \right] \bigg.$$

$$ \times \left[ (2n-3) x-  A_n^2(t) - A_{n-1}^2(t)- {b_n(t)\over a_n(t)}-  {b_{n-1}(t)\over a_{n-1} (t)}\right.$$

$$\bigg.\bigg.+  \left[ x+B_{n-1}(t)-\nu+1-2n\right] [x- B_{n-1}(t)]\bigg]\bigg]  {\partial\over \partial x} [ P_n^\nu(x,t)] $$

$$+ x \bigg[   x \left( n( x - n-\nu ) - {b_n(t)\over a_n(t)} \right) +  \left( B_n(t)-\nu-2n-1\right) \left( A^2_n(t) +  n (x-n-\nu-1) \right) \bigg.$$

$$\bigg. +  \left[  A_n^2(t) +{b_n(t)\over a_n(t)} \right] \left(2n+\nu - B_n(t)\right) \bigg] P_n^\nu(x,t) = 0.$$
Dividing by $x$ and appealing to (3.4), we arrive at (3.6).

\end{proof} 

{\bf Remark 2.} For the limit case $t=0$  we employ (1.5), (2.12), (2.13) to reduce (3.6) to the classical differential equation for Laguerre polynomials.

The following theorem gives the integro-differential-difference  equation  for  orthogonal polynomials $P_n^\nu(x,t)$.  Precisely, it has

{\bf Theorem 5.}  {\it Orthogonal polynomials $P_n^\nu(x,t)$ satisfy the integral-difference equation of the form}

$$  P_n^\nu(x,t)  =   {P_{n-1}^\nu(x,t)\over a_{n-1,0}(t)}  -  a_n(t) \exp\left(  {b_n(t)\over x\ a_n(t)} \right) \int_0^t  \exp\left( -  {b_n(y)\over x\ a_n(y)}  \right) {d\over dy}  \left[ {P_{n-1}^\nu(x,y)\over a_{n-1,0}(y)\ a_n(y)} \right]dy $$

$$+ (-1)^n (n-1)! \ a_n(t) \exp\left(  {1\over x}\left[ {b_n(t)\over a_n(t)} +n(n+\nu)\right]\right) \left[n\ L_n^\nu (x) - \Gamma(1+\nu)  \left( {n! (n+\nu) \over \Gamma(n+\nu)}\right)^{1/2}  L_{n-1}^\nu (x)\right].\eqno(3.9)$$

\begin{proof} Recalling the differential-difference equations  (2.5), (2.28),  the recurrence relation (2.11) and Corollary 2, we rewrite equation (2.22) in the form

$$ {\partial\over \partial t} \left[ {P_n^\nu(x,t) \over a_{n,0}(t)} \right]  + {1\over t} \left[ n  +  {A_n^2(t)  + {b_n(t)\over a_n(t)}\over  B_n(t)-\nu-1-2n } - {1\over x} \left[ A_n^2(t) + {b_n(t)\over a_n(t)} \right]  \right]{P_n^\nu(x,t) \over a_{n,0}(t)}  $$

$$ =  - {1 \over xt} \left[ A_n^2(t) + {b_n(t)\over a_n(t)} \right] {P_{n-1}^\nu(x,t)\over a_{n-1,0}(t)},\quad x \neq 0.\eqno(3.10)$$
Then by virtue of (2.26), (2.37) it implies

$$ {\partial\over \partial t} \left[ {P_n^\nu(x,t) \over a_{n,0}(t)} \right]  + \left[  {a^\prime_{n,0}(t)\over a_{n,0}(t)} -  {a^\prime_{n}(t)\over a_{n}(t)} - {1\over x} {d\over dt}  \left[{b_n(t)\over a_n(t)}\right] \right] {P_n^\nu(x,t) \over a_{n,0}(t)}   =  - {1 \over x} {d\over dt}  \left[{b_n(t)\over a_n(t)}\right]  {P_{n-1}^\nu(x,t)\over a_{n-1,0}(t)}.$$
Solving this first order differential equation in terms of $P_n^\nu$ under initial condition (1.2) and taking into account (2.12), we find

$$ P_n^\nu(x,t)  = -{ a_n(t)\over x}  \exp\left(  {b_n(t)\over x\ a_n(t)} \right) \int_0^t  \exp\left( -  {b_n(y)\over x\ a_n(y)}  \right) {d\over dy}  \left[{b_n(y)\over a_n(y)} \right] {P_{n-1}^\nu(x,y)\over a_{n-1,0}(y)\ a_n(y)} dy $$

$$+ (-1)^n n! \ a_n(t) \exp\left(  {1\over x}\left[ {b_n(t)\over a_n(t)} +n(n+\nu)\right]\right) L_n^\nu (x).\eqno(3.11)$$
Finally, integrating by parts in (3.11), we end up with (3.9).  Theorem 5 is proved.

\end{proof} 

Concerning recurrence relations for the coefficients of the polynomials $P_n^\nu$, we have the following result.

{\bf Theorem 6}. {\it For the orthogonal polynomial $P_n^\nu(x,t) = \sum_{k=0}^n a_{n,k}(t) x^k,\  a_{n,n}(t) \equiv a_n(t),\ a_{n,n-1}(t)= b_n(t)$ its coefficients fullfil the differential-recurrence relations}

$$ a_{n,k}(t)=  a_{n,0}(t) \sum_{m=k}^{n} {  a_{m} (t)  \left[ a_{m}(t) a^\prime_{m,k-1}(t)  -   a_{m,k-1}(t)  a_{m}^\prime(t) \right]  \over \left[ a_{m}(t) b_{m}^\prime(t)  - b_{m}(t) a_{m}^\prime(t) \right] a_{m,0}(t)},\quad k=1,\dots, n.\eqno(3.12)$$

\begin{proof}  By virtue of (2.26) and taking into account that $A_n^2(t)  + {b_n(t)\over a_n(t)} \neq 0,\  t > 0$ (see (2.1), (3.2)), we write equality (3.10) in the form

$$ { a^2_n(t) \over \left[  a_n(t) b_n^\prime(t)  - b_n(t) a_n^\prime(t) \right] a_{n,0}(t) } \left[  {\partial\over \partial t} \left[ P_n^\nu(x,t) \right] - {a_n^\prime(t)\over a_{n}(t)} P_n^\nu(x,t) \right]=    {1 \over x} \left[ {P_n^\nu(x,t)\over a_{n,0}(t)} - {P_{n-1}^\nu(x,t)\over a_{n-1,0}(t)}\right].$$
Hence it yields 

$$   \sum_{m=1}^n {  a_m (t) \over \left[ a_m(t) b_m^\prime(t)  - b_m(t) a_m^\prime(t) \right] a_{m,0}(t)} \sum_{k=1}^{m} \left[ a_m(t) a^\prime_{m,k-1}(t)  -   a_{m,k-1}(t) a_m^\prime(t) \right] x^{k} =    {P_n^\nu(x,t)\over a_{n,0}(t)} -  1.$$
Therefore changing the order of summation on the left-hand side of the latter equality, we arrive at the representation

$$ P_n^\nu(x,t) =   a_{n,0}(t) +   \sum_{k=1}^{n}  x^k  \sum_{m=k}^{n} {  a_{m} (t)  \left[ a_{m}(t) a^\prime_{m,k-1}(t)  -   a_{m,k-1}(t)  a_{m}^\prime(t) \right] a_{n,0}(t) \over \left[ a_{m}(t) b_{m}^\prime(t)  - b_{m}(t) a_{m}^\prime(t) \right] a_{m,0}(t)},$$
which leads to (3.12).

\end{proof}

The three term recurrence relation (1.4) can be written in another form. Indeed, it has

{\bf Theorem 7.} {\it Orthogonal polynomials  $P_n^\nu$ satisfy the following recurrence relation}

$$ \bigg[  A_{n+1}^2(t) + {b_{n+1}(t)\over a_{n+1}(t)} + {x\over 2} \left[B_{n}(t)-\nu-1-2n\right] \bigg] \bigg[ A_{n+1} (t) P_{n+1}^\nu(x,t)+ A_n(t) P_{n-1}^\nu(x,t)\bigg] $$

$$+  \bigg[x t B^\prime_n(t)+ A^2_{n} (t) \left[B_{n-1}(t)-\nu+1-2n\right] - A^2_{n+1}(t) \left[ B_{n+1}(t)-\nu-3-2n\right] \bigg.$$

$$\bigg. - \left(x-  B_n(t)   \right) \left[   A_n^2(t) + {b_n(t)\over a_n(t)}  + {x\over 2} \left[B_n(t)-\nu-1-2n\right] \right]\bigg] P_{n}^\nu(x,t)= 0.\eqno(3.13)$$

\begin{proof}     Differentiating the  three term recurrence relation (1.4), we have

$$x  {\partial\over \partial t} [P_n^\nu(x,t)]  = A^\prime_{n+1} (t) P_{n+1}^\nu(x,t) +B^\prime_n(t) P_n^\nu(x,t) + A^\prime_{n}(t) P_{n-1}^\nu(x,t)$$

$$+   A_{n+1} (t)   {\partial\over \partial t} [P_{n+1}^\nu(x,t)] +B_n(t)  {\partial\over \partial t} [P_n^\nu(x,t)] + A_{n}(t)  {\partial\over \partial t} [P_{n-1}^\nu(x,t)].$$
Then, using (1.4), (2.5) and Corollary 2, after straightforward simplifications the latter equality becomes

$$   \bigg[xt A^\prime_{n+1} (t)+  A_{n +1} (t) \left[   A_{n+1}^2(t) + {b_{n+1}(t)\over a_{n+1}(t)} + {x\over 2} \left[B_{n+1}(t)-\nu-3-2n\right] \right] \bigg]P_{n+1}^\nu(x,t) $$

$$+  \bigg[ xt B^\prime_n(t)+  A^2_{n} (t) \left[ B_{n-1}(t)-\nu+1-2n\right] - A^2_{n+1}(t) \left[ B_{n+1}(t)-\nu-3-2n\right] \bigg.$$

$$\bigg. - \left(x-  B_n(t)   \right) \left[   A_n^2(t) + {b_n(t)\over a_n(t)}  + {x\over 2} \left[B_n(t)-\nu-1-2n\right] \right]\bigg] P_{n}^\nu(x,t)$$

$$\bigg. + \bigg[  xt A^\prime_{n}(t)+  A_{n}(t) \left[   A_{n-1}^2(t) + {b_{n-1}(t)\over a_{n-1}(t)} + \left[B_{n-1}(t) - {x\over 2} \right]\left[ B_{n-1}(t)-\nu+1-2n\right]  \right]\bigg.$$

$$\bigg. -   \left(B_n(t)  - x\right)  A_n(t) \left[ B_n(t)-\nu-1-2n\right]  \bigg]P_{n-1}^\nu(x,t) = 0.$$
It can be rewritten,  appealing to (2.36), and we obtain

$$\bigg[  A_{n+1}^2(t) + {b_{n+1}(t)\over a_{n+1}(t)} + {x\over 2} \left[B_{n}(t)-\nu-1-2n\right] \bigg] A_{n+1} (t) P_{n+1}^\nu(x,t) $$

$$+  \bigg[xB_n(t)+ A^2_{n} (t) \left[x+ B_{n-1}(t)-\nu+1-2n\right] - A^2_{n+1}(t) \left[ x+B_{n+1}(t)-\nu-3-2n\right] \bigg.$$

$$\bigg. - \left(x-  B_n(t)   \right) \left[   A_n^2(t) + {b_n(t)\over a_n(t)}  + {x\over 2} \left[B_n(t)-\nu-1-2n\right] \right]\bigg] P_{n}^\nu(x,t)$$

$$ + \bigg[  A_{n-1}^2(t) + {b_{n-1}(t)\over a_{n-1}(t)} + B_{n-1}(t) \left[ B_{n-1}(t)-\nu+1-2n\right] \bigg.$$

$$\bigg.  -   \left(B_n(t)  - {x\over 2} \right) \left[ B_n(t)-\nu-1-2n\right]  \bigg] A_n(t) P_{n-1}^\nu(x,t) = 0.\eqno(3.14)$$
Hence,  taking into account (2.27), (3.4), equality (3.14) leads to the final form (3.15). 

\end{proof}

{\bf Corollary 8}. {\it The following identity takes place}

$$t B_n(t) B^\prime_n(t)+ A^2_{n} (t) \left[ B_{n-1}(t)-\nu+1-2n\right] - A^2_{n+1}(t) \left[ B_{n+1}(t)-\nu-3-2n\right]  = 0.\eqno(3.16)$$

\begin{proof} The proof is immediate, letting $x=0$ in (3.14) and recalling (2.27), (3.4).

\end{proof}

{\bf Corollary 9}.  {\it Coefficients $B_n$ obey the following integral-recurrence relation}

$$ B_n(t)+ B_{n-1}(t) - 2n-\nu +1 =  \exp\left( \int_0^t {B_{n-1}(y)-B_{n}(y)+2\over y} dy\right) $$

$$\times \left( 2\int_0^t  \exp\left( - \int_0^y {B_{n-1}(u)-B_{n}(u)+2\over u} du\right)  \left[ B_{n-1}(y) - 2n-\nu +1\right] {dy\over y}+ 2n+\nu+1\right).\eqno(3.17)$$

\begin{proof} In fact, writing (3.3) in the form

$${d\over dt}  \bigg[ B_n(t)+ B_{n-1}(t) - 2n-\nu +1\bigg] + {1\over t} \left[ B_{n-1}(t)-B_{n}(t)+2\right] \left[ B_{n-1}(t)+ B_{n}(t) - 2n-\nu +1\right] $$

$$- {2\over t} \left[ B_{n-1}(t) - 2n-\nu +1\right] = 0,$$
we solve a simple Cauchy problem   for the first order differential equation to obtain (3.17).

\end{proof}

In the sequel, let us consider polynomial coefficients (3.12) $a_{n,k}(t)$ as functions of $\nu$ as well, i.e.  $a_{n,k}\equiv a_{n,k}^\nu$.   Then, returning to the formula (1.8) for the moments, we  represent  orthogonal polynomials $P_n^\nu$ in terms of the Hankel determinant

$$ P^\nu_{n} (x,t)= {1\over [G^\nu_{n-1} (t)G^\nu_n (t)]^{1/2} } \begin{vmatrix}  

\rho_{\nu+1}(t)  & \rho_{\nu+2}(t)&  \dots&  \dots&   \rho_{\nu+n+1}(t) \\

\rho_{\nu+2}(t) &  \dots &   \dots&    \dots&    \rho_{\nu+n+2}(t) \\
 
\dots&  \dots &   \dots&   \dots&   \dots \\
 
 \vdots  &   \ddots  &  \ddots & \ddots&   \vdots\\
  
 \rho_{\nu+n}(t) &   \dots&  \dots&   \dots&    \rho_{\nu+2n}(t)\\

1 &   x & \dots&    \dots&    x^n\\

 \end{vmatrix},\eqno(3.18)$$
where

$$G^\nu_{n}(t) =   \begin{vmatrix}  

\rho_{\nu+1}(t)  & \rho_{\nu+2}(t)&  \dots&  \dots&   \rho_{\nu+n+1}(t) \\

\rho_{\nu+2}(t) &  \dots &   \dots&    \dots&    \rho_{\nu+n+2}(t) \\
 
\dots&  \dots &   \dots&   \dots&   \dots \\
 
 \vdots  &   \ddots  &  \ddots & \ddots&   \vdots\\
  
 \rho_{\nu+n+1}(t) &   \dots&  \dots&   \dots&    \rho_{\nu+2n+1}(t)\\

 \end{vmatrix}.\eqno(3.19)$$
 Hence, employing the Laplace theorem to the last row, we find from (3.18), (3.19)
 
 $$G^\nu_{n}(t) = (-1)^n \rho_{\nu+n+1}(t) G^{\nu+1}_{n-1}(t) +  \rho_{\nu+2n+1}(t) G^{\nu}_{n-1}(t)+ (-1)^{n+1} \sum_{k=2}^{n}  (-1)^k \rho_{\nu+n+k}(t) G^\nu_{n,k}(t)$$
 
 $$=   [G^\nu_{n-1} (t)G^\nu_n (t)]^{1/2}  \sum_{k=0}^{n}   \rho_{\nu+n+k+1}(t) a^\nu_{n,k}(t),\eqno(3.20)$$
  where 
  
  $$G^\nu_{n,k}(t) =  \begin{vmatrix}  

\rho_{\nu+1}(t)  & \rho_{\nu+2}(t)&  \dots&  \rho_{\nu+k-1}(t)&  \rho_{\nu+k+1}(t) &  \dots& \rho_{\nu+n+1}(t) \\

\rho_{\nu+2}(t) &  \rho_{\nu+3}(t) &   \dots&  \rho_{\nu+k}(t)&  \rho_{\nu+k+2}(t) &  \dots& \rho_{\nu+n+2}(t) \\
 
\dots&  \dots &   \dots&   \dots&   \dots &\dots  &\dots  \\
 
 \vdots  &   \vdots  &  \vdots & \ddots&   \vdots &\vdots &\vdots\\
  
 \rho_{\nu+n}(t) &   \dots&  \dots&   \rho_{\nu+n+k-2}(t)&  \rho_{\nu+n+k}(t) &  \dots& \rho_{\nu+2n}(t)\\

 \end{vmatrix}.\eqno(3.21)$$
 In particular,  we easily find from (3.18), (3.20) the expression for the free coefficient $a_{n,0}^\nu(t)$ in terms of the determinant by the formula

$$ a_{n,0}^\nu(t)  =  {(-1)^n  G^\nu_{n,1}(t) \over \left[G^\nu_{n-1} (t)G^\nu_n (t)\right]^{1/2} }.\eqno(3.22)$$
Moreover, since the leading term $a_n^\nu$ (cf. \cite{Bateman}, Vol. II) has the representation

$$ \left|a_n^\nu(t) \right| =  \bigg[ {G^{\nu}_{n-1}(t) \over G^{\nu}_{n}(t)}\bigg]^{1/2},\eqno(3.23)$$
we find from (3.22), (3.23) the equality 

$$ {G^{\nu+1}_{n-1}(t) \over G^{\nu}_{n}(t)} = (-1)^n  a_n^\nu(t)\ a_{n,0}^\nu(t).\eqno(3.24)$$
Therefore it yields

$$ (- 1)^n  a_n^\nu(t)\ a_{n,0}^\nu(t) = {G^{\nu+1}_{n-1}(t) \over G^{\nu+1}_{n}(t)}  {G^{\nu+1}_{n}(t) \over G^{\nu}_{n}(t)} $$

$$= (- 1)^{n+1} \left[ a_n^{\nu+1}(t)\right]^2 \ a_{n+1,0}^\nu(t) a_{n+1}^\nu(t)  {G^{\nu}_{n+1}(t) \over G^{\nu}_{n}(t)}$$

$$= (- 1)^{n+1} \left[ {a_n^{\nu+1}(t)\over a_{n+1}^\nu(t)}\right]^2 \ a_{n+1,0}^\nu(t) a_{n+1}^\nu(t).$$
Hence, recalling (2.11), we derive the identities 

$$  \left[ {a_n^{\nu+1}(t)\over a_{n+1}^\nu(t)}\right]^2 = - A_{n+1}^\nu(t) {a_{n,0}^\nu(t) \over a_{n+1,0}^\nu(t)} = \frac {\left[A^\nu_{n+1}(t)\right]^2+ {b^\nu_{n+1}(t)\over a^\nu_{n+1}(t)}}{2n+\nu+3- B^\nu_{n+1}(t) },$$

$$ \left[ a_n^{\nu+1}(t)\right]^2 = \frac {\left[a^\nu_{n}(t)\right]^2+ a^\nu_{n+1}(t) b^\nu_{n+1}(t)}{2n+\nu+3- B^\nu_{n+1}(t) },\eqno(3.25)$$

$$\prod_{k=0}^n \left[ a_k^{\nu+1}(t)\right]^2 =   {(-1)^{n+1} a_{n+1}^\nu(t)\over a_{n+1,0}^\nu(t)} \prod_{k=0}^n \left[ a_k^{\nu}(t)\right]^2.\eqno(3.26)$$
Moreover,  (3.20), (3.23) imply

$$ \sum_{k=0}^{n}   \rho_{\nu+n+k+1}(t) a^\nu_{n,k}(t) = {1\over a_n^\nu(t)},\eqno(3.27)$$
and this agrees with (1.20). 

{\bf Corollary 10.} {\it Coefficients $B^\nu_n(t), A^\nu_{n+1}(t)$ $(1.5)$ of the three term recurrence relation $(1.4)$ are equal, correspondingly,}

$$ B^\nu_n(t)= 2n+\nu+1 + t \bigg[ {{G^{\nu}}^\prime_{n-1}(t) \over G^{\nu}_{n-1}(t) } - {{G^{\nu}}^\prime_{n}(t) \over G^{\nu}_{n}(t)}\bigg],\eqno(3.28)$$

$$A_{n+1}^\nu(t)= {\left[G^{\nu}_{n-1}(t) G^{\nu}_{n+1}(t) \right]^{1/2}\over G^{\nu}_{n}(t)}.\eqno(3.29)$$

\begin{proof} Indeed, appealing to (2.16) and (3.23),  we have

$$ {G^{\nu}_{n-1}(t) \over G^{\nu}_{n}(t)} =  {1 \over n! \Gamma(n+\nu+1)}\ \exp\left(\int_0^t {B_n(y)-\nu-1-2n\over y} dy\right).\eqno(3.30)$$
Therefore,

$$G^{\nu}_{n}(t) =  n! \Gamma(n+\nu+1)\  G^{\nu}_{n-1}(t) \exp\left(\int_0^t {\nu+1+2n - B_n(y)\over y} dy\right),$$
and, solving this recurrence relation owing to (1.5) and (3.19), we easily derive 

$$G^{\nu}_{n}(t) =  \prod_{k=0}^n  k! \Gamma(k+\nu+1)  \exp\left(\int_0^t \bigg[ {b_{n+1}^\nu(y)\over a_{n+1}^\nu(y)} + (n+1)(n+\nu+1) \bigg] {dy \over y} \right).$$
Consequently, 

$${G^{\nu}}^\prime_{n}(t) \equiv {dG^{\nu}_{n}\over dt} = {G^{\nu}_{n}(t) \over t}   \bigg[ {b_{n+1}^\nu(t)\over a_{n+1}^\nu(t)} + (n+1)(n+\nu+1) \bigg].\eqno(3.31)$$
Hence equality (3.28) follows immediately from (1.5). Formula (3.29) is a direct  consequence of (1.5), (3.23) or (2.18), (3.30).
\end{proof} 

{\bf Corollary 11.} {\it The following equality holds}

$$t^2 {d\over dt}\bigg[ {{G^{\nu}}^\prime_{n}\over  G^{\nu}_{n}}\bigg]= A_{n+1}^2 -  (n+1)(n+1+\nu).\eqno(3.32)$$

\begin{proof} The proof is immediate, involving differentiation in (3.31) and employing (2.26).

\end{proof}

On the other hand, the key identity (1.16) for the moments and properties of the determinants allow to treat (3.19) when $n = 2,3,\dots,$ as follows

$$  t G^\nu_{n}(t) =   \begin{vmatrix}  

t \rho_{\nu+1}(t)  & \rho_{\nu+2}(t)&  \dots&  \dots&   \rho_{\nu+n+1}(t) \\

t \rho_{\nu+2}(t) &  \dots &   \dots&    \dots&    \rho_{\nu+n+2}(t) \\
 
\dots&  \dots &   \dots&   \dots&   \dots \\
 
 \vdots  &   \ddots  &  \ddots & \ddots&   \vdots\\
  
 t \rho_{\nu+n+1}(t) &   \dots&  \dots&   \dots&    \rho_{\nu+2n+1}(t)\\

 \end{vmatrix}$$

$$  =  - \begin{vmatrix}  

(\nu+2) \rho_{\nu+2}(t)  & \rho_{\nu+2}(t)&  \dots&  \dots&   \rho_{\nu+n+1}(t) \\

(\nu+3) \rho_{\nu+3}(t) &  \rho_{\nu+3}(t) &   \dots&    \dots&    \rho_{\nu+n+2}(t) \\
 
\dots&  \dots &   \dots&   \dots&   \dots \\
 
 \vdots  &   \vdots  &  \ddots & \ddots&   \vdots\\
  
  (\nu+n+2)\rho_{\nu+n+2}(t) &   \rho_{\nu+n+2}(t)&  \dots&   \dots&    \rho_{\nu+2n+1}(t)\\

 \end{vmatrix}$$

$$  =  - \begin{vmatrix}  

0 & \rho_{\nu+2}(t)&  \dots&  \dots&   \rho_{\nu+n+1}(t) \\

\rho_{\nu+3}(t) &  \rho_{\nu+3}(t)  &   \dots&    \dots&    \rho_{\nu+n+2}(t) \\
 
2 \rho_{\nu+4}(t) &  \rho_{\nu+4}(t)  &   \dots&   \dots&   \dots \\
 
 \vdots  &   \vdots  &  \ddots & \ddots&   \vdots\\
  
  n \rho_{\nu+n+2}(t) &   \rho_{\nu+n+2}(t) &  \dots&   \dots&    \rho_{\nu+2n+1}(t)\\

 \end{vmatrix}.$$
Continuing this process,  we find via Laplace's theorem and (1.16)

$$(-1)^{n+1} t^{n-1} G^\nu_n(t) = \begin{vmatrix}  

0 & 0&  \dots& 0& \rho_{\nu+n}(t)&   \rho_{\nu+n+1}(t) \\

\rho_{\nu+3}(t) &  \rho_{\nu+4}(t)  &   \dots&     \rho_{\nu+n+1}(t)&    \rho_{\nu+n+1}(t) &  \rho_{\nu+n+2}(t)\\
 
2 \rho_{\nu+4}(t) &   2 \rho_{\nu+5}(t) &   \dots&   2 \rho_{\nu+n+2}(t)&    \rho_{\nu+n+2}(t)&    \rho_{\nu+n+3}(t) \\
 
 \vdots  &   \vdots  &  \ddots & \ddots&   \vdots & \vdots\\
  
  n \rho_{\nu+n+2}(t) &   n\rho_{\nu+n+3}(t) &  \dots&    n \rho_{\nu+2n}(t)&  \rho_{\nu+2n}(t)&    \rho_{\nu+2n+1}(t)\\

 \end{vmatrix}.\eqno(3.33)$$
  On the other hand, using continuously (1.16) from the last column in (3.19),  we arrive in the same manner at the equality

$$ G^\nu_n(t) = \begin{vmatrix}  

 \rho_{\nu+1}(t)&   \rho_{\nu+2}(t) &0 & 0&  \dots& 0 \\

 \rho_{\nu+2}(t) &  \rho_{\nu+3}(t)& \rho_{\nu+3}(t) &  \rho_{\nu+4}(t)  &   \dots&     \rho_{\nu+n+1}(t)\\
 
   \rho_{\nu+3}(t)&    \rho_{\nu+4}(t) & 2 \rho_{\nu+4}(t) &   2 \rho_{\nu+5}(t) &   \dots&   2 \rho_{\nu+n+2}(t) \\
 
 \vdots  &   \vdots  &  \vdots & \vdots&   \vdots &   \vdots \\
  
    \rho_{\nu+n+1}(t)&    \rho_{\nu+n+2}(t)&  n \rho_{\nu+n+2}(t) &   n\rho_{\nu+n+3}(t) &  \dots&    n \rho_{\nu+2n}(t)\\

 \end{vmatrix}.\eqno(3.34)$$
  Let us consider the following double sequence of determinants
 
 $$H_{i,j}^{\nu,n}(t)= \begin{vmatrix}  

 \rho_{\nu+i}(t)&   \rho_{\nu+j}(t) &0 & 0&  \dots& 0 \\

 \rho_{\nu+i+1}(t) &  \rho_{\nu+j+1}(t)& \rho_{\nu+3}(t) &  \rho_{\nu+4}(t)  &   \dots&     \rho_{\nu+n+1}(t)\\
 
   \rho_{\nu+i+2}(t)&    \rho_{\nu+j+2}(t) & 2 \rho_{\nu+4}(t) &   2 \rho_{\nu+5}(t) &   \dots&   2 \rho_{\nu+n+2}(t) \\
 
 \vdots  &   \vdots  &  \vdots & \vdots&   \vdots &   \vdots \\
  
    \rho_{\nu+n+i}(t)&    \rho_{\nu+n+j}(t)&  n \rho_{\nu+n+2}(t) &   n\rho_{\nu+n+3}(t) &  \dots&    n \rho_{\nu+2n}(t)\\

 \end{vmatrix},\eqno(3.35)$$
where $(i,j) \in \mathbb{Z}^2$. Hence it is not difficult to observe from (3.33), (3.34), (3.35) that $H_{i,j}^{\nu,n}(t)= - H_{j,i}^{\nu,n}(t)$ and 

$$H_{j,j}^{\nu,n}(t)=0, \ j \in \mathbb{Z}, \quad   H_{j,j+1}^{\nu,n}(t) =  (-1)^{j+1} t^{j-1} G^\nu_n(t),\quad j=1,\dots,n.\eqno(3.36)$$
Generally, we find via (1.16)

$$   H_{i,j}^{\nu,n}(t)= (\nu+i-1)  H_{i-1,j}^{\nu,n}(t) + t H_{i-2,j}^{\nu,n}(t)$$

$$+ \begin{vmatrix}  

 0& \rho_{\nu+j}(t)&  0 &0 &  \dots& 0 \\

 \rho_{\nu+i}(t) &  \rho_{\nu+j+1}(t)& \rho_{\nu+3}(t) &  \rho_{\nu+4}(t)  &   \dots&     \rho_{\nu+n+1}(t)\\
 
  2 \rho_{\nu+i+1}(t)&    \rho_{\nu+j+2}(t) & 2 \rho_{\nu+4}(t) &   2 \rho_{\nu+5}(t) &   \dots&   2 \rho_{\nu+n+2}(t) \\
 
 \vdots  &   \vdots  &  \vdots & \vdots&   \vdots &   \vdots \\
  
    n \rho_{\nu+n+i-1}(t)&    \rho_{\nu+n+j}(t)&  n \rho_{\nu+n+2}(t) &   n\rho_{\nu+n+3}(t) &  \dots&    n \rho_{\nu+2n}(t)\\

 \end{vmatrix}$$
 
$$= (\nu+i-1)  H_{i-1,j}^{\nu,n}(t) + t H_{i-2,j}^{\nu,n}(t)$$

$$- n! \rho_{\nu+j}(t) \begin{vmatrix}  

 \rho_{\nu+i}(t) &  \rho_{\nu+3}(t) &  \rho_{\nu+4}(t)  &   \dots&     \rho_{\nu+n+1}(t)\\
 
   \rho_{\nu+i+1}(t)&   \rho_{\nu+4}(t) &    \rho_{\nu+5}(t) &   \dots&    \rho_{\nu+n+2}(t) \\
 
 \vdots  &   \vdots  &  \vdots & \vdots&      \vdots \\
  
  \rho_{\nu+n+i-1}(t)&    \rho_{\nu+n+2}(t) &   \rho_{\nu+n+3}(t) &  \dots&     \rho_{\nu+2n}(t)\\

 \end{vmatrix}.$$
Hence,

$$   H_{i,j}^{\nu,n}(t)= (\nu+i-1)  H_{i-1,j}^{\nu,n}(t) + t H_{i-2,j}^{\nu,n}(t)$$

$$- n! \rho_{\nu+j}(t) \begin{vmatrix}  

 \rho_{\nu+i}(t) &  \rho_{\nu+3}(t) &  \rho_{\nu+4}(t)  &   \dots&     \rho_{\nu+n+1}(t)\\
 
   \rho_{\nu+i+1}(t)&   \rho_{\nu+4}(t) &    \rho_{\nu+5}(t) &   \dots&    \rho_{\nu+n+2}(t) \\
 
 \vdots  &   \vdots  &  \vdots & \vdots&      \vdots \\
  
  \rho_{\nu+n+i-1}(t)&    \rho_{\nu+n+2}(t) &   \rho_{\nu+n+3}(t) &  \dots&     \rho_{\nu+2n}(t)\\

 \end{vmatrix}.\eqno(3.37)$$
Analogously, it has 

$$   H_{i,j}^{\nu,n}(t)= (\nu+j-1)  H_{i,j-1}^{\nu,n}(t) + t H_{i,j-2}^{\nu,n}(t)$$

$$+ n! \rho_{\nu+i}(t) \begin{vmatrix}  

 \rho_{\nu+j}(t) &  \rho_{\nu+3}(t) &  \rho_{\nu+4}(t)  &   \dots&     \rho_{\nu+n+1}(t)\\
 
   \rho_{\nu+j+1}(t)&   \rho_{\nu+4}(t) &    \rho_{\nu+5}(t) &   \dots&    \rho_{\nu+n+2}(t) \\
 
 \vdots  &   \vdots  &  \vdots & \vdots&      \vdots \\
  
  \rho_{\nu+n+j-1}(t)&    \rho_{\nu+n+2}(t) &   \rho_{\nu+n+3}(t) &  \dots&     \rho_{\nu+2n}(t)\\

 \end{vmatrix},\eqno(3.38)$$
Thus we find

$$(\nu+i-1)  H_{i-1,j}^{\nu,n}(t) + t H_{i-2,j}^{\nu,n}(t)$$

$$- n! \rho_{\nu+j}(t) \begin{vmatrix}  

 \rho_{\nu+i}(t) &  \rho_{\nu+3}(t) &  \rho_{\nu+4}(t)  &   \dots&     \rho_{\nu+n+1}(t)\\
 
   \rho_{\nu+i+1}(t)&   \rho_{\nu+4}(t) &    \rho_{\nu+5}(t) &   \dots&    \rho_{\nu+n+2}(t) \\
 
 \vdots  &   \vdots  &  \vdots & \vdots&      \vdots \\
  
  \rho_{\nu+n+i-1}(t)&    \rho_{\nu+n+2}(t) &   \rho_{\nu+n+3}(t) &  \dots&     \rho_{\nu+2n}(t)\\

 \end{vmatrix}$$

$$= (\nu+j-1)  H_{i,j-1}^{\nu,n}(t) + t H_{i,j-2}^{\nu,n}(t)$$

$$+ n! \rho_{\nu+i}(t) \begin{vmatrix}  

 \rho_{\nu+j}(t) &  \rho_{\nu+3}(t) &  \rho_{\nu+4}(t)  &   \dots&     \rho_{\nu+n+1}(t)\\
 
   \rho_{\nu+j+1}(t)&   \rho_{\nu+4}(t) &    \rho_{\nu+5}(t) &   \dots&    \rho_{\nu+n+2}(t) \\
 
 \vdots  &   \vdots  &  \vdots & \vdots&      \vdots \\
  
  \rho_{\nu+n+j-1}(t)&    \rho_{\nu+n+2}(t) &   \rho_{\nu+n+3}(t) &  \dots&     \rho_{\nu+2n}(t)\\

 \end{vmatrix}.\eqno(3.39)$$
As  an immediate consequence, when $(i, j) \in [ 3,\dots, n+1] \times [ 3,\dots, n+1]$ the following recurrence relation holds

$$(\nu+i-1)  H_{i-1,j}^{\nu,n}(t) + (\nu+j-1)  H_{j-1,i}^{\nu,n}(t) + t\bigg[ H_{i-2,j}^{\nu,n}(t) + H_{j-2,i}^{\nu,n}(t)\bigg] = 0.\eqno(3.40)$$
Recalling  (3.36), (3.37), it   yields

$$H_{n+2,n+1}^{\nu,n}(t) = (-1)^{n+1} \bigg[  t^{n-1}  G^{\nu}_{n}(t) - n! \rho_{\nu+n+1}(t) G^{\nu+2}_{n-1}(t)\bigg].\eqno(3.41) $$
On the other hand,  differentiating determinant (3.19) by virtue of (1.15) and applying the same process as above to the derivative,  we deduce the equalities

  $$t {d G^\nu_{n}\over dt}  = -   \begin{vmatrix}  

t \rho_{\nu}(t)  & \rho_{\nu+2}(t)&  \dots&  \dots&   \rho_{\nu+n+1}(t) \\

t\rho_{\nu+1}(t) &  \dots &   \dots&    \dots&    \rho_{\nu+n+2}(t) \\
 
\dots&  \dots &   \dots&   \dots&   \dots \\
 
 \vdots  &   \ddots  &  \ddots & \ddots&   \vdots\\
  
t \rho_{\nu+n}(t) &   \dots&  \dots&   \dots&    \rho_{\nu+2n+1}(t)\\

 \end{vmatrix}$$
 
$$= (\nu+1) G_n^\nu(t) +  \begin{vmatrix}  

0 & \rho_{\nu+2}(t)&  \dots&  \dots&   \rho_{\nu+n+1}(t) \\

\rho_{\nu+2}(t) &  \rho_{\nu+3}(t)  &   \dots&    \dots&    \rho_{\nu+n+2}(t) \\
 
2 \rho_{\nu+3}(t) &  \rho_{\nu+4}(t)  &   \dots&   \dots&   \dots \\
 
 \vdots  &   \ddots  &  \ddots & \ddots&   \vdots\\
  
  n \rho_{\nu+n+1}(t) &   \rho_{\nu+n+2}(t) &  \dots&   \dots&    \rho_{\nu+2n+1}(t)\\

 \end{vmatrix}$$

$$ =  (\nu+1) G_n^\nu(t) +  \begin{vmatrix}  

0 & \rho_{\nu+2}(t)&  \rho_{\nu+3}(t) & 0& \dots& 0&  0 \\

\rho_{\nu+2}(t) &  \rho_{\nu+3}(t)  &    \rho_{\nu+4}(t)&  \rho_{\nu+4}(t)&\dots&   \rho_{\nu+n}(t) &    \rho_{\nu+n+1}(t) \\
 
2 \rho_{\nu+3}(t) &  \rho_{\nu+4}(t)  &    \rho_{\nu+5}(t)&  2\rho_{\nu+5}(t)& \dots&  2 \rho_{\nu+n}(t) &   2 \rho_{\nu+n+1}(t ) \\
 
 \vdots  &   \ddots  &  \ddots & \ddots&   \vdots&\vdots &\vdots\\
  
  n \rho_{\nu+n+1}(t) &   \rho_{\nu+n+2}(t) & \rho_{\nu+n+3}(t) & n \rho_{\nu+n+3}(t)& \dots&   n \rho_{\nu+n-1}(t)&   n \rho_{\nu+2n}(t)\\

 \end{vmatrix}$$

 $$ =  (\nu+1) G_n^\nu(t) +  \begin{vmatrix}  

0 & \rho_{\nu+2}(t)& 0 & 0& \dots& 0&  0 \\

\rho_{\nu+2}(t) &  \rho_{\nu+3}(t)  &    \rho_{\nu+3}(t)&  \rho_{\nu+4}(t)&\dots&   \rho_{\nu+n}(t) &    \rho_{\nu+n+1}(t) \\
 
2 \rho_{\nu+3}(t) &  \rho_{\nu+4}(t)  &   2 \rho_{\nu+4}(t)&  2\rho_{\nu+5}(t)& \dots&  2 \rho_{\nu+n}(t) &   2 \rho_{\nu+n+1}(t ) \\
 
 \vdots  &   \ddots  &  \ddots & \ddots&   \vdots&\vdots &\vdots\\
  
  n \rho_{\nu+n+1}(t) &   \rho_{\nu+n+2}(t) & n \rho_{\nu+n+2}(t) & n \rho_{\nu+n+3}(t)& \dots&   n \rho_{\nu+n-1}(t)&   n \rho_{\nu+2n}(t)\\

 \end{vmatrix}$$

$$+  t  \begin{vmatrix}  

0 & \rho_{\nu+2}(t)&  \rho_{\nu+1}(t) & 0& \dots& 0&  0 \\

\rho_{\nu+2}(t) &  \rho_{\nu+3}(t)  &    \rho_{\nu+2}(t)&  \rho_{\nu+4}(t)&\dots&   \rho_{\nu+n}(t) &    \rho_{\nu+n+1}(t) \\
 
2 \rho_{\nu+3}(t) &  \rho_{\nu+4}(t)  &    \rho_{\nu+3}(t)&  2\rho_{\nu+5}(t)& \dots&  2 \rho_{\nu+n}(t) &   2 \rho_{\nu+n+1}(t ) \\
 
 \vdots  &   \ddots  &  \ddots & \ddots&   \vdots&\vdots &\vdots\\
  
  n \rho_{\nu+n+1}(t) &   \rho_{\nu+n+2}(t) & \rho_{\nu+n+1}(t) & n \rho_{\nu+n+3}(t)& \dots&   n \rho_{\nu+n-1}(t)&   n \rho_{\nu+2n}(t)\\

 \end{vmatrix}.$$
Consequently, it gives by virtue of (3.19), (1.16) 

  $$t {d G^\nu_{n}\over dt} =  (\nu+1) G_n^\nu(t) - n!  \rho_{\nu+2}(t) G_{n-1}^{\nu+1}(t) $$

$$+  t^2  \begin{vmatrix}  

 \rho_{\nu}(t)&  \rho_{\nu+1}(t)&0  & 0& \dots& 0&  0 \\

 \rho_{\nu+1}(t)  &    \rho_{\nu+2}(t)& \rho_{\nu+2}(t)& \rho_{\nu+4}(t)&\dots&   \rho_{\nu+n}(t) &    \rho_{\nu+n+1}(t) \\
 
  \rho_{\nu+2}(t)  &    \rho_{\nu+3}(t)&  2 \rho_{\nu+3}(t) & 2\rho_{\nu+5}(t)& \dots&  2 \rho_{\nu+n}(t) &   2 \rho_{\nu+n+1}(t ) \\
 
 \vdots  &   \vdots  &  \vdots & \vdots&   \vdots&\vdots &\vdots\\
  
   \rho_{\nu+n}(t) & \rho_{\nu+n+1}(t) &n \rho_{\nu+n+1}(t)  & n \rho_{\nu+n+3}(t)& \dots&   n \rho_{\nu+n-1}(t)&   n \rho_{\nu+2n}(t)\\

 \end{vmatrix}.\eqno(3.42)$$
 But from (3.35) we observe

 $${d\over dt} \bigg[ H_{0,1}^{\nu,n}(t) \bigg] = - \begin{vmatrix}  

 \rho_{\nu-1}(t)&   \rho_{\nu+1}(t) &0 & 0&  \dots& 0 \\

 \rho_{\nu}(t) &  \rho_{\nu+2}(t)& \rho_{\nu+3}(t) &  \rho_{\nu+4}(t)  &   \dots&     \rho_{\nu+n+1}(t)\\
 
   \rho_{\nu+1}(t)&    \rho_{\nu+3}(t) & 2 \rho_{\nu+4}(t) &   2 \rho_{\nu+5}(t) &   \dots&   2 \rho_{\nu+n+2}(t) \\
 
 \vdots  &   \vdots  &  \vdots & \vdots&   \vdots &   \vdots \\
  
    \rho_{\nu+n-1}(t)&    \rho_{\nu+n+1}(t)&  n \rho_{\nu+n+2}(t) &   n\rho_{\nu+n+3}(t) &  \dots&    n \rho_{\nu+2n}(t)\\

 \end{vmatrix}$$
 
 $$- \begin{vmatrix}  

 \rho_{\nu}(t)&   \rho_{\nu+1}(t) &0 & 0&  \dots& 0 \\

 \rho_{\nu+1}(t) &  \rho_{\nu+2}(t)& \rho_{\nu+2}(t) &  \rho_{\nu+4}(t)  &   \dots&     \rho_{\nu+n+1}(t)\\
 
   \rho_{\nu+2}(t)&    \rho_{\nu+3}(t) & 2 \rho_{\nu+3}(t) &   2 \rho_{\nu+5}(t) &   \dots&   2 \rho_{\nu+n+2}(t) \\
 
 \vdots  &   \vdots  &  \vdots & \vdots&   \vdots &   \vdots \\
  
    \rho_{\nu+n}(t)&    \rho_{\nu+n+1}(t)&  n \rho_{\nu+n+1}(t) &   n\rho_{\nu+n+3}(t) &  \dots&    n \rho_{\nu+2n}(t)\\

 \end{vmatrix}.$$
 Hence we derive from (3.42) the identity 
 
 $$t {d G^\nu_{n}\over dt} -  (\nu+1) G_n^\nu(t) + n!  \rho_{\nu+2}(t) G_{n-1}^{\nu+1}(t) + t^2\bigg[ {d\over dt} \bigg[ H_{0,1}^{\nu,n}(t) \bigg] + H_{-1,1}^{\nu,n}(t) \bigg] = 0.\eqno(3.43)$$ 
{\bf Remark 3}.   By using these features an interesting open question is to obtain a differential-difference recurrence relation for determinants  (3.19).

\section{Rodrigues-type formula. Generating function}

Let us expand the function $e^{-t/x} P^\nu_n(x,t)$ in terms of the Laguerre polynomials $L_n^\nu(x)$.  It gives

$$ e^{-t/x} P^\nu_n(x,t) = \sum_{k=0}^\infty d_{n,k}^\nu(t) L_k^\nu(x),\eqno(4.1)$$
where

$$ d_{n,k}^\nu(t) = {k!\over \Gamma(k+\nu+1)} \int_0^\infty  e^{-x -t/x} P^\nu_n(x,t) L_k^\nu(x) x^\nu dx.\eqno(4.2)$$
But from (1.18) we have $d_{n,k}^\nu(t) = 0, \ k =0,\dots, n-1$.  The uniform estimate with respect to $k$ for coefficients $d_{n,k}^\nu(t)$ is given by the following lemma.

{\bf Lemma 2}. {\it Let $ t >0\ \nu > -1$.  Coefficients $d_{n,k}^\nu(t),\  n,k  \in \mathbb{N}_0$ satisfy the upper bound of the form

$$\left|d_{n,k}^\nu(t)\right| \le  {k!\  h_n^\nu(t) \over \Gamma(k+\nu+1)},\eqno(4.3)$$
where

$$h_n^\nu(t)=   2^{\nu-1/2} \int_0^t Q_n^\nu(t-y) \rho^{1/2}_{2\nu+1} (2 y) {dy\over \sqrt y}\eqno(4.4)$$
and }

$$Q_n^\nu(x) = \sum_{m=0}^n \left| a_{n,m}^\nu(t)\right| { x^m\over m!}.\eqno(4.5)$$

\begin{proof} Writing polynomial $P^\nu_n$ in the explicit form, equality (4.2) becomes

$$ d_{n,k}^\nu(t) = {k!\over \Gamma(k+\nu+1)} \sum_{m=0}^n a_{n,m}^\nu(t) \int_0^\infty  e^{-x -t/x} L_k^\nu(x) x^{\nu+m} dx.\eqno(4.6)$$
Recalling (1.11), the latter integral can be represented as follows

$$ \int_0^\infty  e^{-x -t/x} L_k^\nu(x) x^{\nu+m} dx = (-1)^{k+m+1} {d^{k-m-1}\over dt^{k-m-1}} \int_0^\infty x^{\nu+k -1} e^{-x - t/x}  L_k^\nu(x) dx$$

$$= {(-1)^{m+1}\over k!} {d^{k-m-1}\over dt^{k-m-1}} \bigg[  t^k \rho_\nu(t) \bigg]= {(-1)^{m+1}\over k!\ m!} \int_0^t (t-y)^m {d^{k}\over dy^{k}} \bigg[  y^k \rho_\nu(y) \bigg] dy,\eqno(4.7)$$
where we mean (cf. (1.12))

$${d^{-q} f\over dt^{-q}} \equiv \left(I_+^q f\right)(t)=  {1\over (q-1)!} \int_0^t (t-y)^{q-1} f(y) dy,\quad q \in \mathbb{N}_0.$$
Taking into account (4.6), coefficients $d_{n,k}^\nu(t)$ can be written in the operator form

$$d_{n,k}^\nu(t) = -  {k!\over \Gamma(k+\nu+1)} P_n^\nu\bigg( - I_+, t\bigg) \bigg\{ {d^{k}\over dt^{k}} \bigg[  t^k \rho_\nu(t) \bigg]\bigg\}.\eqno(4.8)$$
Meanwhile, via (1.11) and the Rodrigues formula for Laguerre polynomials it has

$$ {d^{k}\over dy^{k}} \bigg[  y^k \rho_\nu(y) \bigg] = k! \int_0^\infty x^{\nu -1} e^{-x - y/x}  L_k^\nu(x) dx =  \int_0^\infty x^{-1} e^{- y/x} {d^k\over dx^k} \left[ e^{-x} x^{\nu+k} \right] dx.$$
Integrating by parts with the use of Entry 1.1.3.2 on p. 4 in \cite{Bry}, we get 

$$  {d^{k}\over dy^{k}} \bigg[  y^k \rho_\nu(y) \bigg]  =  (-1)^k \int_0^\infty  {d^k\over dx^k} \left[ x^{-1} e^{- y/x} \right]  e^{-x} x^{\nu+k} dx$$

$$=  k!  \int_0^\infty   e^{-x- y/x} x^{\nu-1} L_k\left({y\over x}\right) dx,$$
where $L_k$ are  Laguerre polynomials of the index zero.   Then owing to (1.7), Schwarz's inequality and orthogonality of Laguerre polynomials, we find the estimate

$$\left|  {d^{k}\over dy^{k}} \bigg[  y^k \rho_\nu(y) \bigg] \right| \le  k!  \left(\int_0^\infty  e^{- yx- 2/x} x^{-2(\nu+1)} dx\right)^{1/2} \left(\int_0^\infty   e^{- yx}  \left[L_k\left(y x\right)\right]^2  dx\right)^{1/2} $$ 

$$= {2^{\nu-1/2} k!\over \sqrt y}  \rho^{1/2}_{2\nu+1} (2 y).$$
Therefore, returning to (4.6), (4.7), we derive

$$\left| d_{n,k}^\nu(t)\right| \le  { 2^{\nu-1/2}  k!\over \Gamma(k+\nu+1)} \int_0^t \sum_{m=0}^n \left| a_{n,m}^\nu(t)\right| { (t-y)^m\over m!}  \rho^{1/2}_{2\nu+1} (2 y) {dy\over \sqrt y},$$
which implies (4.3) and completes the proof of Lemma 2.

\end{proof} 

{\bf Corollary 12.}  {\it Under the condition $\nu > 3/2$ the Laguerre series $(4.1)$ converges absolutely and uniformly on closed intervals of $\mathbb{R}_+$.}

\begin{proof}  In fact,  Stirling's asymptotic formula for gamma-function \cite{Bateman}, Vol. I yields

$$ {k!\over \Gamma(k+\nu+1)} = O\left( k^{-\nu} \right),\quad k \to \infty.$$
Since Laguerre polynomials $L_k^\nu(x)$ behave as $O\left( k^{\nu/2- 1/4}\right),\ k \to \infty$ uniformly on closed  intervals $x \in [\alpha, \beta]$ of $\mathbb{R}_+$, the absolute and uniform convergence of the series (4.1) is guaranteed under the assumption $\nu > 3/2$.

\end{proof}

Further,  writing (4.1)  as follows

$$ e^{-t/x} P^\nu_n(x,t) = \sum_{k=n}^\infty d_{n,k}^\nu(t) L_k^\nu(x) = x^{-\nu} e^x \sum_{k=n}^\infty {d_{n,k}^\nu(t)\over k!} {d^k\over dx^k } \left[x^{k+\nu} e^{-x} \right]$$

$$= x^{-\nu} e^x \sum_{k=0}^\infty {d_{n,k+n}^\nu(t)\over (k+n)!}  {d^{k+n}\over dx^{k+n} } \left[x^{k+n+\nu} e^{-x} \right]$$

$$=  x^{-\nu} e^x   \sum_{k=0}^\infty {d_{n,k+n}^\nu(t)\ k! \over (k+n)!} {d^{n}\over dx^{n} } \bigg[ e^{-x} x^{\nu+n}  L_k^{n+\nu} (x)\bigg],$$
the problem arises to interchange the order of differentiation and summation.   It can be done by virtue of the absolute and uniform convergence of the series 

$$\sum_{k=0}^\infty {d_{n,k+n}^\nu(t)\over (k+n)!}  {d^{k+j}\over dx^{k+j} } \left[x^{k+n+\nu} e^{-x} \right]$$

$$=  x^{-\nu-n+j} e^x \sum_{k=0}^\infty {d_{n,k+n}^\nu(t)(k+j)! \over (k+n)!} L_{k+j}^{n+\nu-j}(x),\quad j=1,2,\dots, n,\ n \in \mathbb{N}\eqno(4.9) $$
on closed intervals of $\mathbb{R}_+$. Indeed, this is an immediate consequence of the bound (4.3) and Corollary 12.  Thus we arrive at the following Rodrigues-type formula for polynomials  $P^\nu_n$

$$   P^\nu_n(x,t) =  x^{-\nu} e^{x+t/x} {d^{n}\over dx^{n} }  \bigg[ e^{-x} x^{\nu+n}  \sum_{k=0}^\infty {d_{n,k+n}^\nu(t)\ k! \over (k+n)!} \  L_k^{n+\nu} (x)\bigg].\eqno(4.10)$$
Recurrence relations for coefficients $d_{n,k}^\nu(t)$ are given by 

{\bf Theorem 8}. {\it Coefficients $d_{n,k}^\nu(t)$ obey  recurrence relations of the form}

$$ \left(2k+\nu+1 - (k+1) B_n^\nu(t) \right) d_{n,k}^\nu(t) -  \left(k+1\right) A_{n+1}^\nu(t) d_{n+1,k}^\nu(t)- ( k+\nu+1) d_{n,k+1}^\nu(t)$$

$$ - \left(k+1\right)  A_{n}^\nu(t) d_{n-1,k}^\nu(t) - k\ d_{n,k-1}^\nu(t) = 0,\eqno(4.11)$$

$$ (2k+n+1+\nu) A_{n+1}^\nu(t) d_{n+1,k}^\nu(t) -  (k+1+\nu) A_{n+1}^\nu(t) d_{n+1,k+1}^\nu(t) $$

$$+ \bigg[ (1+n)  B_n^\nu(t) - (k+\nu+2)(2k+\nu+1) -t +\left[ A_n^\nu(t)\right]^2- A_n^\nu(t) \left[ B^\nu_n(t)-\nu-1-2n\right] \bigg] d_{n,k}^\nu(t) $$

$$+ (k+\nu+1)(k+\nu+2- B^\nu_n(t) ) d_{n,k+1}^\nu(t) -  (k+\nu+1) A_n^\nu(t) d_{n-1,k+1}^\nu(t)  $$

$$+\bigg[ \left(2k+n+\nu+1+ B^\nu_{n-1}(t) \right) A_n^\nu(t) - \bigg[  [A^\nu_n(t)]^2 + {b^\nu_n(t)\over a^\nu_n(t)} \bigg] \bigg] d_{n-1,k}^\nu(t) $$

$$+ k(k+\nu+2) d_{n,k-1}^\nu(t)+  A_n^\nu(t) A_{n-1}^\nu(t) d_{n-2,k}^\nu(t) = 0.\eqno(4.12)$$

\begin{proof}  Indeed, relation (4.11) follows immediately from the three recurrence relations (1.3), (1.4)  for Laguerre polynomials and polynomials $P_n^\nu$, respectively.  In order to prove (4.12), we integrate by parts in (4.2) to obtain

$$d_{n,k}^\nu(t)=  {k!\over t \Gamma(k+\nu+1)} \bigg[ \int_0^\infty  e^{-x -t/x} P^\nu_n(x,t) L_k^\nu(x) x^{\nu+2} dx\bigg.$$

$$ - (\nu+2) \int_0^\infty  e^{-x -t/x} P^\nu_n(x,t) L_k^\nu(x) x^{\nu+1} dx + \int_0^\infty  e^{-x -t/x} P^\nu_n(x,t) L_{k-1}^{\nu+1}(x) x^{\nu+2} dx$$

$$\bigg. - \int_0^\infty  e^{-x -t/x} {\partial\over \partial x} \left[ P^\nu_n(x,t)\right]  L_k^\nu(x) x^{\nu+2} dx\bigg].\eqno(4.13)$$
But using the known relation for Laguerre polynomials \cite{Sze}

$$x L_{k-1}^{\nu+1}(x) = (k+\nu) L_{k-1}^\nu(x) - k L_k^\nu(x),$$
equality (4.13) becomes

$$d_{n,k}^\nu(t)=  {k!\over t \Gamma(k+\nu+1)} \bigg[ \int_0^\infty  e^{-x -t/x} P^\nu_n(x,t) L_k^\nu(x) x^{\nu+2} dx\bigg.$$

$$ - (\nu+2+k) \int_0^\infty  e^{-x -t/x} P^\nu_n(x,t) L_k^\nu(x) x^{\nu+1} dx + (k+\nu) \int_0^\infty  e^{-x -t/x} P^\nu_n(x,t) L_{k-1}^{\nu}(x) x^{\nu+1} dx$$

$$\bigg. - \int_0^\infty  e^{-x -t/x} {\partial\over \partial x} \left[ P^\nu_n(x,t)\right]  L_k^\nu(x) x^{\nu+2} dx\bigg].\eqno(4.14)$$
Then, employing again recurrence relations (1.3), (1.4) and differential-difference equation (2.5), we derive

$$  {k!\over t \Gamma(k+\nu+1)} \int_0^\infty  e^{-x -t/x} P^\nu_n(x,t) L_k^\nu(x) x^{\nu+2} dx$$

$$=  {k!\over t \Gamma(k+\nu+1)}  \int_0^\infty  e^{-x -t/x} \left[ A_{n+1}^\nu(t) P^\nu_{n+1}(x,t) + B_{n}^\nu(t) P^\nu_{n}(x,t) + A_{n}^\nu(t) P^\nu_{n-1}(x,t)\right] $$

$$\times \left[ (2k+1+\nu) L_k^\nu(x)  - (k+1)L_{k+1}^\nu(x) - (k+\nu)L_{k-1}^\nu(x)\right] x^\nu dx $$

$$= {1\over t} \bigg[  (2k+1+\nu) \bigg[ A_{n+1}^\nu(t) d_{n+1,k}^\nu(t) + B_{n}^\nu(t) d_{n,k}^\nu(t) + A_{n}^\nu(t) d_{n-1,k}^\nu(t)\bigg] \bigg.$$

$$ -  (k+1+\nu) \bigg[ A_{n+1}^\nu(t) d_{n+1,k+1}^\nu(t) + B_{n}^\nu(t) d_{n,k+1}^\nu(t) + A_{n}^\nu(t) d_{n-1,k+1}^\nu(t)\bigg]$$

$$\bigg. -  k \bigg[ A_{n+1}^\nu(t) d_{n+1,k-1}^\nu(t) + B_{n}^\nu(t) d_{n,k-1}^\nu(t) + A_{n}^\nu(t) d_{n-1,k-1}^\nu(t)\bigg] \bigg],$$

$$ {k! (k+\nu+2) \over t \Gamma(k+\nu+1)} \int_0^\infty  e^{-x -t/x} P^\nu_n(x,t) L_k^\nu(x) x^{\nu+1} dx $$

$$= {k+\nu+2 \over t} \bigg[  (2k+1+\nu) d_{n,k}^\nu(t) - (k+1+\nu) d_{n,k+1}^\nu(t) - k d_{n,k-1}^\nu(t)\bigg],$$

$$  {k!\over t \Gamma(k+\nu)} \int_0^\infty  e^{-x -t/x} P^\nu_n(x,t) L_{k-1}^{\nu}(x) x^{\nu+1} dx $$

$$= {k\over t} \bigg[ A_{n+1}^\nu(t) d_{n+1,k-1}^{\nu} (t) +  B_{n}^\nu(t) d_{n,k-1}^{\nu} (t) + A_{n}^\nu(t)  d_{n-1,k-1}^{\nu} (t)\bigg],$$

$$ {k!\over t \Gamma(k+\nu+1)} \int_0^\infty  e^{-x -t/x} {\partial\over \partial x} \left[ P^\nu_n(x,t)\right]  L_k^\nu(x) x^{\nu+2} dx $$

$$= {k!\over t \Gamma(k+\nu+1)} \int_0^\infty  e^{-x -t/x}  L_k^\nu(x) x^{\nu} \bigg[ \left[ nx - \left[A^\nu_n(t)\right]^2  - {b^\nu_n(t)\over a^\nu_n(t)} \right] P_n^\nu(x,t)\bigg.$$

$$\bigg. + A^\nu_n(t) \left[ x+B^\nu_n(t)-\nu-1-2n\right] P_{n-1}^\nu(x,t) \bigg] dx$$

$$=   {A^\nu_n(t)\over t}  \left[ B^\nu_n(t)-\nu-1-2n\right]  d_{n,k}^\nu(t) - {1\over t} \bigg[ \left[A^\nu_n(t)\right]^2 + {b^\nu_n(t)\over a^\nu_n(t)} \bigg] d_{n-1,k}^\nu(t)$$

$$+ {n\over t}  \bigg[ A_{n+1}^\nu(t) d_{n+1,k}^\nu(t) + B_{n}^\nu(t) d_{n,k}^\nu(t) + A_{n}^\nu(t) d_{n-1,k}^\nu(t)\bigg]$$

$$+ {A^\nu_n(t) \over t} \bigg[ A_{n}^\nu(t) d_{n,k}^\nu(t) + B_{n-1}^\nu(t) d_{n-1,k}^\nu(t) + A_{n-1}^\nu(t) d_{n-2,k}^\nu(t)\bigg].$$
Hence,  substituting these values into (4.13), we get after simplification the desired relation (4.12).

\end{proof}

Finally,  defining   as usual the generating function $G(x,w,t)$ in terms of the series

$$G(x,w,t)= \sum_{n=0}^\infty  P^\nu_n(x,t) {w^n\over n!},\quad x >0, \ w \in \mathbb{C},\eqno(4.15)$$
it can be written due to (4.10) in the form 

$$G(x,w,t)=    x^{-\nu} e^{x+t/x} \sum_{n=0}^\infty  {w^n\over n!}   {d^{n}\over dx^{n} }  \bigg[ e^{-x} x^{\nu+n}  \sum_{k=0}^\infty {d_{n,k+n}^\nu(t)\ k! \over (k+n)!} \  L_k^{n+\nu} (x)\bigg].\eqno(4.16)$$
The convergence of the series (4.15) is guaranteed at least in $L_2\left(\mathbb{R}_+;  e^{-x-t/x} x^\nu dx\right)$.  Indeed, by virtue of the Minkowski inequality it has 

$$\bigg(\int_0^\infty e^{-x-t/x} \bigg|   \sum_{n=N}^\infty  P^\nu_n(x,t)  {w^n\over n!} \bigg|^2  x^\nu dx \bigg)^{1/2} \le  \sum_{n=N}^\infty {|w|^n\over n!} \bigg(\int_0^\infty e^{-x-t/x} \left[P^\nu_n(x,t)\right]^2  x^\nu dx \bigg)^{1/2} $$

$$=   \sum_{n=N}^\infty {|w|^n\over n!} \to 0,\quad N \to \infty.$$

\bibliographystyle{amsplain}

\end{document}